\documentclass[11pt]{article}
\usepackage{amsmath,amssymb,graphicx,url}

\newcommand{\RR}{\mathbb{R}}

\newcommand{\QQ}{\mathbb{Q}}

\newcommand{\CC}{\mathbb{C}}
\newtheorem{Theorem}{Theorem}
\newtheorem{Proposition}[Theorem]{Proposition}

\begin{document}
\title{Elementary fractal geometry. \\ 
Networks and carpets involving irrational rotations} 
\author{Christoph Bandt and Dmitry Mekhontsev}
\maketitle

\begin{abstract}
Self-similar sets with open set condition, the linear objects of fractal geometry, have been considered mainly for crystallographic data. Here we introduce new symmetry classes in the plane, based on rotation by irrational angles. Examples without characteristic directions, with strong connectedness and small complexity were found in a computer-assisted search. They are surprising since the rotations are given by rational matrices, and the proof of the open set condition usually requires integer data. 
We develop a classification of self-similar sets by symmetry class and algebraic numbers. Examples are given for various quadratic number fields. . 
\end{abstract}
 
\section{Introduction}\label{intro}
\begin{figure}[h!] 
\begin{center}
\includegraphics[width=0.32\textwidth]{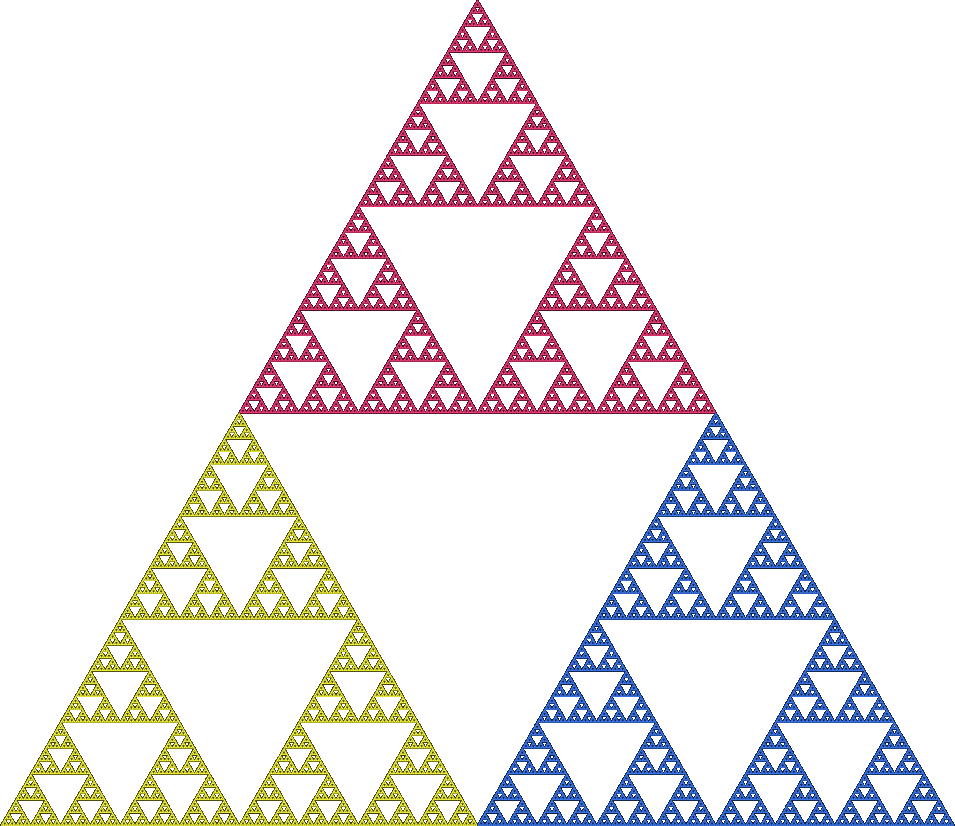}\quad
\includegraphics[width=0.29\textwidth]{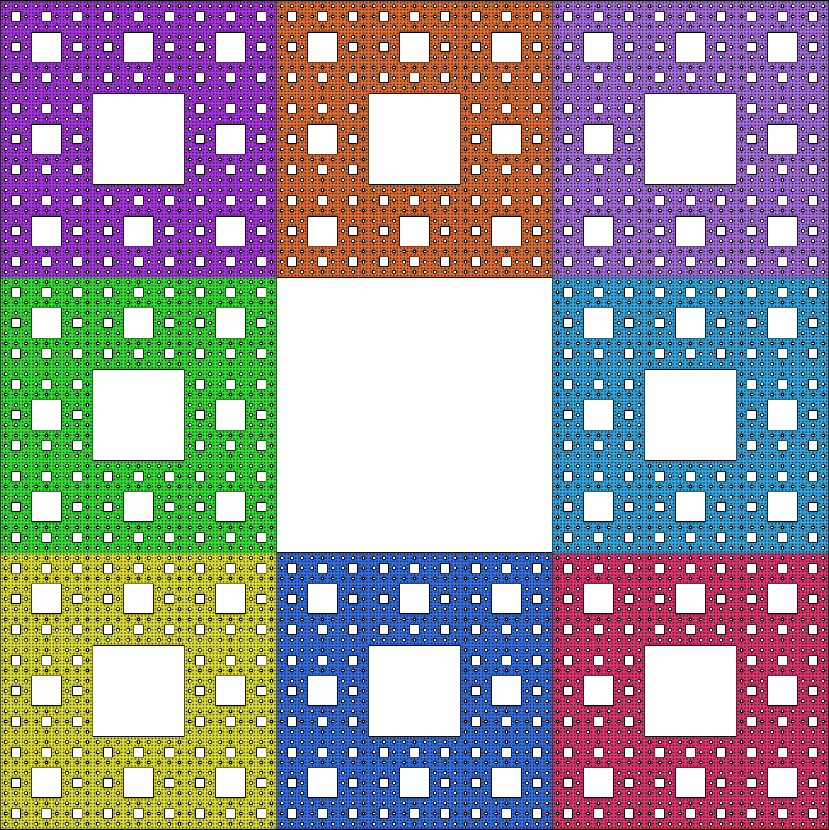}\quad
\includegraphics[width=0.32\textwidth]{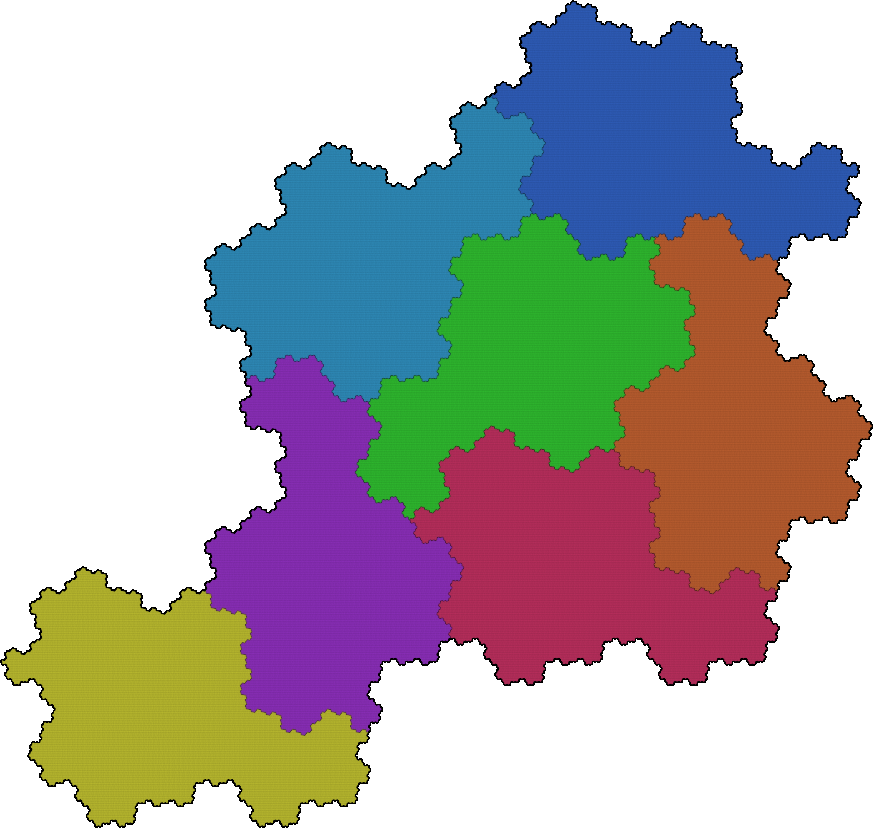}
\end{center}
\caption{Classical fractal shapes. The Sierpi\'nski triangle, the Sierpi\'nski carpet and a variation of the Gosper snowflake tile are standard examples of self-similar sets with open set condition.}\label{fi1}
\end{figure}

\subparagraph{Self-similar sets and open set condition.}
A self-similar set is a nonempty compact subset $A$ of $\RR^d$ which is the union of shrinked copies of itself. This is expressed by Hutchinson's equation
\begin{equation} A=\bigcup_{k=1}^m f_k(A)  
\label{hut}\end{equation}
Here $F=\{ f_1,...,f_m\}$ is a finite set of contractive similarity mappings, often called an iterated function system, abbreviated IFS. Recall that a contractive similarity map from Euclidean $\RR^d$ to itself fulfils $|f(x)-f(y)|=r_f|x-y|$ where the constant  $r_f<1$ is called the factor of $f.$ To keep things simple, we assume that all maps in $F$ have the same factor $r.$  For a given IFS  $F=\{ f_1,...,f_m\}$ there is a unique self-similar set $A$ which is often called the attractor of $F.$ See \cite{Bar,BP,Fal} for details. Figure \ref{fi1} shows standard examples with factor $\frac12, \frac13 ,$ and $\sqrt{1/7},$ respectively.  

The idea of equation \eqref{hut} is that $A$ subdivides into $m^2$ pieces $f_{k_1}f_{k_2}(A)$ of second level, and into $m^n$ pieces $f_w(A)$ of level $n,$ where $w=k_1...k_n$ runs through all words of length $n$ with letters from the alphabet $K=\{ 1,...,m\} .$ Thus we have a homogeneous structure of little pieces, and we can define a uniform measure on $A$ by assigning each pieces of level $n$ the value $m^{-n}.$

In 1946, Moran \cite{Mo} constructed the uniform measure on $A$ as normalized Hausdorff measure of dimension $\alpha$ where $m r^\alpha =1.$ 
He needed a condition, the so-called open set condition, or OSC for short: 
\begin{equation} 
\mbox{There is an open } U \mbox{ so that the image sets } f_1(U),...,f_m(U)    
\mbox{ are disjoint subsets of } U.
\label{osc}\end{equation}
In Figure \ref{fi1}, $U$ can be taken as an open triangle, a square, and the interior of the set $A,$ respectively.
It has turned out that without open set condition, the compact set $A$ does not have a nice local structure. The number of sister pieces in the vicinity of a piece of level $n$ will tend to infinity with $n$ if OSC does not hold \cite{BG,Sch,BP}. However, self-similar sets play a similar part in fractal geometry as lines in ordinary geometry.  The concept of self-similarity is that after sufficient magnification, the view of $A$ will repeat and should not become infinitely dense. That is why OSC is required. In the literature we find mainly examples like Figure \ref{fi1} where the open set can be easily constructed.

\subparagraph{The aim of this paper.}
Sierpi\'nski constructed his triangle and carpet more than 100 years ago as topological spaces with curious properties \cite{Man}. After 1980, physicists took them as models of porous materials, and mathematicians developed an analysis of heat equation, Brownian motion, and eigenvalues of the Laplace operator on just these spaces. See the books by Kigami and Strichartz \cite{Kig,Str} for an introduction, and the literature cited there.
However, both examples have special properties which are rarely met in nature:
\begin{itemize}
\item They have a few characteristic directions.
\item They contain line segments.
\item Their holes have small perimeter compared to their area.
\end{itemize}
Moreover, Sierpi\'nski's triangle and the related p.c.f. fractals \cite{Kig,Str} have cutpoints (see Section \ref{ngapp}) and resemble networks more than porous materials. The vector space of harmonic functions on such spaces is finite-dimensional. The topology of the carpet is definitely more realistic and interesting, but also hard to study. 

Nowadays, complex geometric structures are studied in nearly every active area of science: cell biology, the brain, soil and roots, foam, clouds, dust, nanostructures etc. Therefore we think that fractal geometry should diversify its models. In this paper we show that this is possible, even within the framework of equation \eqref{hut} with equal contraction factors in the plane. This is a beginning. We hope to study more general cases in a subsequent paper.

Here we construct carpets based on non-crystallographic IFS data which 
\begin{itemize}
\item have no characteristic directions and an isotropic type of symmetry,
\item contain no line segments, and
\item have holes with large perimeter and rather complicated shape.
\end{itemize}

In statistical physics, much more complicated random fractals were introduced by Sheffield and Werner \cite{WS}.  Their conformal loop ensemble is a probability space of fractal shapes which is invariant not only under particular similitudes, but under arbitrary conformal maps. Instances of the ensemble are hard to imagine and still harder to visualize. Concrete examples in this paper can be seen as a step from Figure \ref{fi1} to such abstract models.  

We shall focus on carpets with small complexity, for which there is a chance to develop fractal analysis. In the next paragraphs, we sketch an approach to complexity of self-similar sets. We define the neighbor graph of an IFS, the automaton which generates the topology of $A.$ The number of states of the automaton is taken as complexity.

\subparagraph{Algebraic OSC.}
There is an algebraic equivalent for the open set condition, formulated explicitly in terms of the IFS $F$ \cite{BG,BP}.
\begin{equation} 
F\mbox{ generates a free semigroup $F^*$, and $id$ is an isolated point in } {F^*}^{-1}F^*  .
\label{nei}\end{equation}
The first part of this condition says that two $n^{\rm th}$ level pieces $f_w(A)$ and $f_v(A)$ will not coincide, for arbitrary $n.$ This part is not essential for our concept of self-similarity, and a weak separation condition WSC was defined by assuming only the second part of condition \eqref{nei} \cite{Zer,NW,LN07,DE11}.
However, in that case we would have to control pieces which appear once, twice, or $n$ times, and we want to keep things simple.  The second part of \eqref{nei} says that the maps $h=f_w^{-1}f_v$ stay away uniformly from the identity map, so that pieces $f_w(A)$ and $f_v(A)$ do not come arbitrary near to the same position, for any level $n.$ Like OSC, this is an accurate formulation of the vague statement that ``the overlap of two different pieces is not too large''.

Things become simpler when we go to a discrete setting. We assume that there is a matrix $M$ associated with our iterated function system $F$ on $\RR^d$ so that each $f_k$ in $F$ has the form
\begin{equation} 
f_k(x) = M^{-1} s_k x + v_k  \mbox{ with a vector }v_k \mbox{ and an orthogonal matrix } s_k .
\label{map}\end{equation}
We can now assume that $v$ is an integer vector and $M, s_k$ are integer matrices which commute with each other. 
Then we have to check only finitely many maps $h=f_w^{-1}f_v,$ and the second part of condition \eqref{nei} will follow from the first one \cite{BG,BP,Zer}.

Equation \eqref{map} is a strong assumption, however. To get similarity maps, $M$ must fulfil $M\cdot M'= \frac1r\cdot I$ where $I$ denotes the unity matrix. And even if we assume that $M$ commutes only with the whole group generated by the $s_k,$ as in Theorem 2 of \cite{Ba5},
there are only a few choices of integer maps $s_k.$ On the positive side, however, we have the fact that all calculations are done in integer arithmetics and thus the check of the OSC is accurate. It was implemented in the package IFStile \cite{M}
and has led to thousands of new examples with surprising properties. A few of them, tightly related to the Sierpi\'nski triangle in Figure \ref{fi1}, were discussed in \cite{Sierrel}.  

\begin{figure}[h!] 
\begin{center}
\includegraphics[width=0.475\textwidth]{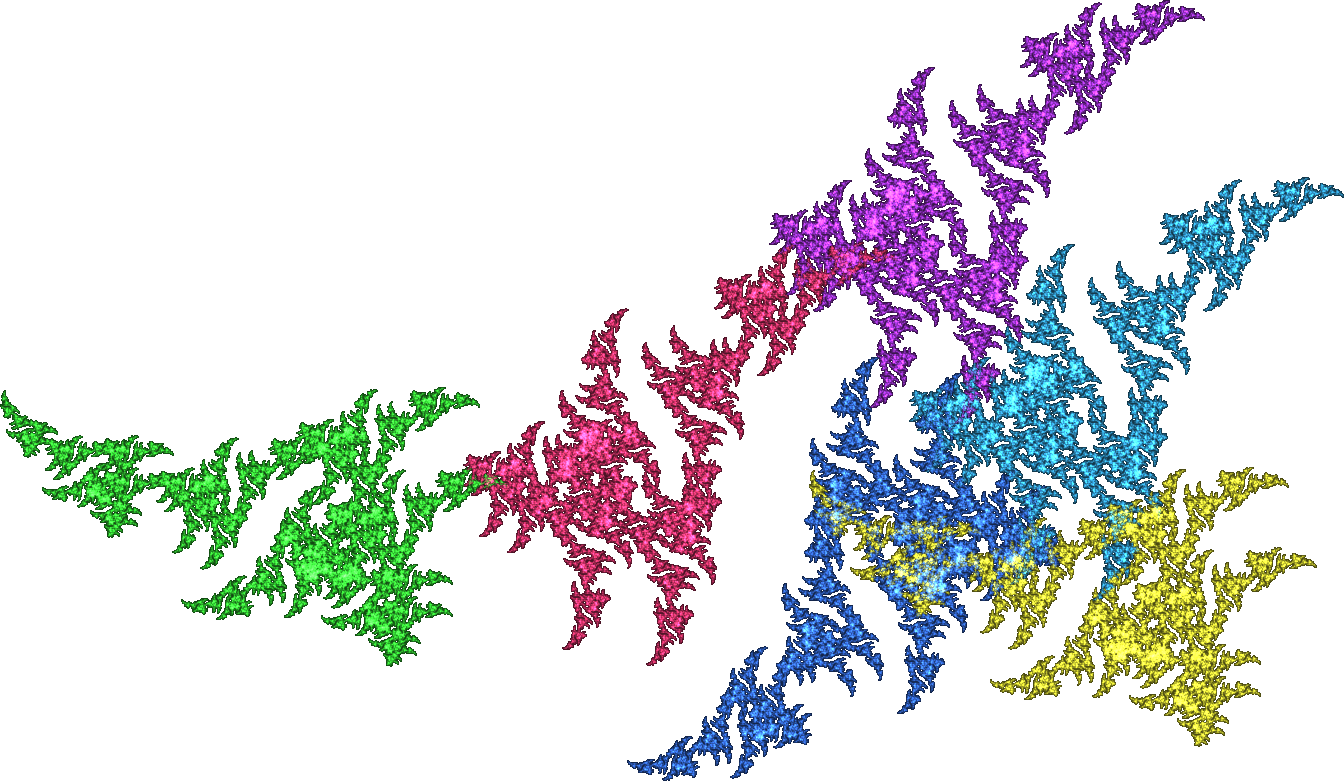} \quad
\includegraphics[width=0.475\textwidth]{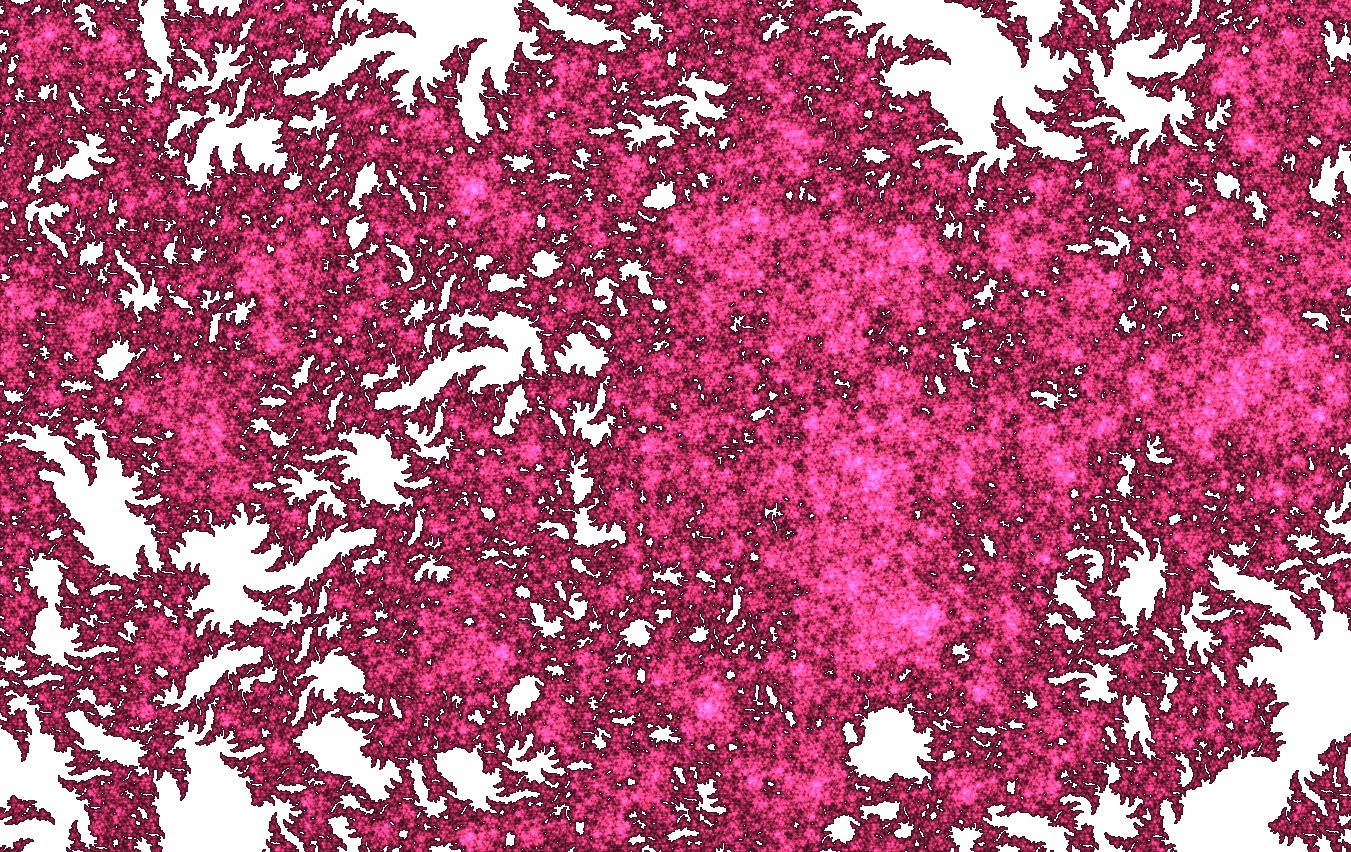}
\end{center}
\caption{An example with six pieces and considerable overlap for which the OSC is still true. In such cases, several magnifications are  necessary to reveal the local structure. The data for this IFS come from a crystallographic group with $60^o$ rotations.}\label{fi2}
\end{figure}  

Another example is given in Figure \ref{fi2}. Since an open set in such examples consists of infinitely many components (cf. \cite{Sierrel}), the local structure of such examples cannot be seen from a global picture of the set $A.$ Usually we need various magnifications of the set in order to understand the local structure. The number of pictures which we need can be seen as a measure of complexity of the self-similar set.  

\subparagraph{Neighbor complexity.}
As a rigorous concept of complexity, we take the number of neighbor maps. Consider the map $h=f_w^{-1}f_v$ where $w,v\in K^n$ denote different words of length $n$ from the alphabet $K.$ It is an isometry. We call $h$ a \emph{neighbor map} or \emph{neighbor type} if the corresponding pieces of $A$ intersect: $f_w(A)\cap f_v(A)\not= \emptyset .$  The number of neighbor types will be taken as the complexity of the IFS $F.$ If the number is finite, we say that $F$ is of finite type. As a consequence, OSC is fulfilled: a special open set was constructed in \cite{BHR}.

The neighbor map  $h=f_w^{-1}f_v$  represents the relative position of the intersecting pieces $f_w(A)$ and $f_v(A).$ 
However, $h$ does not map the pieces into each other. This is done by the map  $g=f_vf_w^{-1}.$ The neighbor map $h$ always maps $A$ to a potential neighbor $h(A)$ in some `superpiece'. Constructions with larger and larger pieces are familiar in the theory of self-similar tilings, see for instance \cite{GS,Se95,Baake2013}. 
The maps $h$ and $g$ are conjugate:  $h=f_v^{-1}gf_v .$ But $g$ depends on the size of the pieces while $h$ is standardized and does not depend on the level $n.$

Suppose the Sierpi\'nski carpet in Figure \ref{fi1} is realized in the complex plane, in a square with vertices $0, 1, 1+i, i. $  Then the neighbor maps will be translations $h(x)=x+v$ with $v\in\{ \pm 1,\pm i\}$ for neighbors with a common edge. For neighbors with a common vertex, we have the translation vectors $v\in \{ \pm (1+i), \pm (1-i)\} .$ So the Sierpi\'nski carpet has 8 neighbor types. The square as a self-similar set with $m=k^2$ pieces forming a $k\times k$ checkerboard pattern has the same neighbor types.  Thus to distinguish different $k$ we need other measures of complexity. The Sierpi\'nski triangle has 6 neighbor types, two for each of its vertices. The tile in Figure \ref{fi1} has 11 types and the example in Figure \ref{fi2} has 52.

With the concept of neighbor type, we provide a quantitative version of the open set condition. In fact we are not interested in the question ``OSC or not OSC''.  We rather want to find examples of small complexity which we can understand.  At present, an IFS with 1000 neighbor types can hardly be distinguished from an IFS without OSC. 

Given an IFS \eqref{map} with integer vectors and matrices, the question whether there are at most $N$ neighbor types is decidable in finite time for every $N.$   One of the authors has written the program IFStile  \cite{M}
which decides this question within milliseconds for $N\le 500.$  The algorithm, discussed in Section \ref{ng}, constructs the neighbor graph of $F,$ an automaton which describes geometric and topological properties of the set $A.$

However, when we assume integer matrices, we are in the setting of crystallographic groups. In the plane, symmetries must be rotations by multiples of 30 degrees, or reflections, the axes of which differ by angles of $k\cdot 15$ degrees. Besides Figure \ref{fi1}, this includes Figure \ref{fi2} where we have no reflection, only rotation by $60^o$ and $180^o.$ Nevertheless, this setting is rather restrictive.

\begin{figure}[h!] 
\begin{center}
\includegraphics[width=0.41\textwidth]{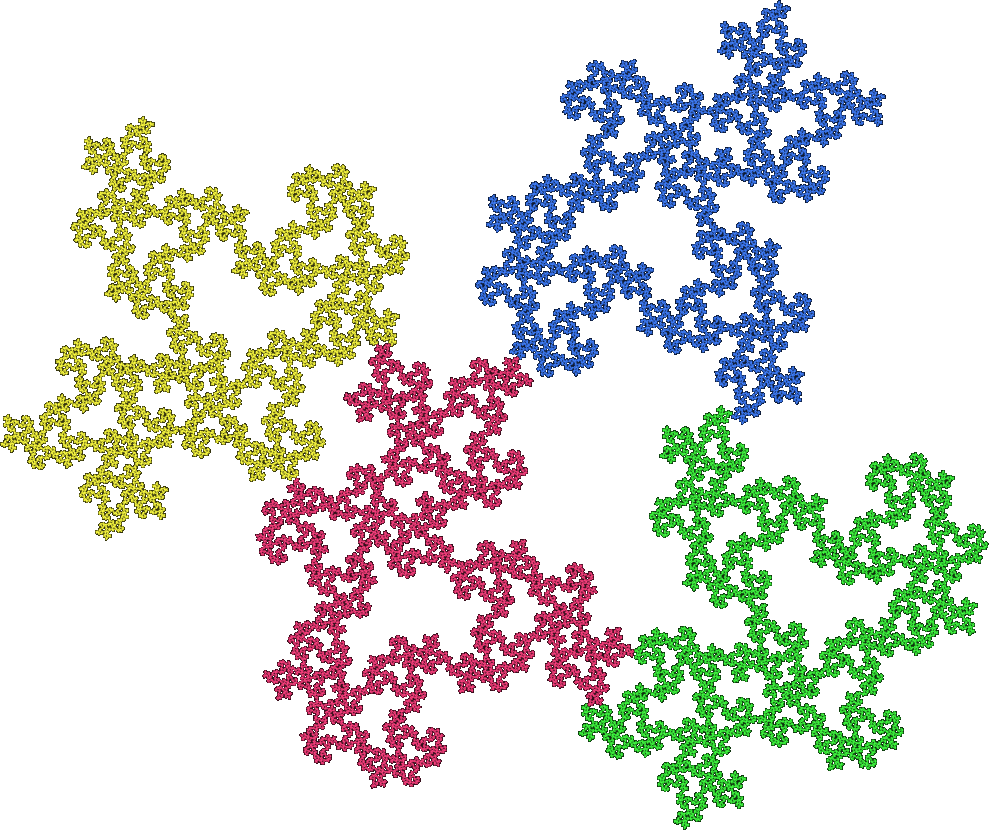} \quad
\includegraphics[width=0.54\textwidth]{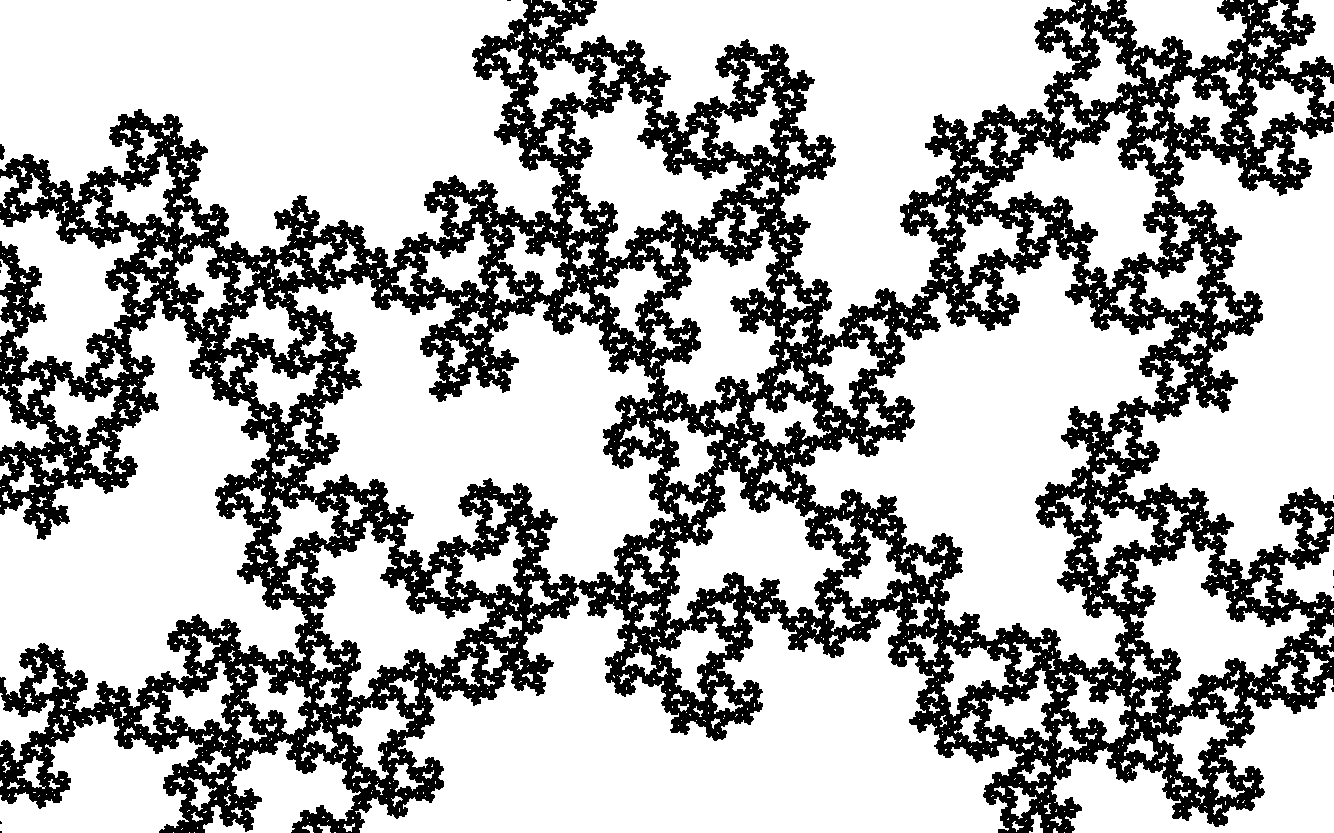}
\end{center}
\caption{A carpet generated from an irrational rotation.  It has only 16 neighbor types, much less than Figure \ref{fi2}. There are no characteristic directions. See the text in Section \ref{pyth}.}\label{fi6}
\end{figure}  

\subparagraph{Approaches to non-crystallographic patterns.}
There are various ways to generalize crystallographic finite type systems. One way is projection of integer data from a higher-dimensional space, known from quasiperiodic tilings used for modelling quasicrystals \cite{Se95, Baake2013}.
This is implemented in IFStile. 
Another possibility is to replace the equation \eqref{hut} by a system of equations for different types of sets, which is called a graph-directed construction \cite{MW}.
In their study of self-similar tilings, Thurston, Kenyon, and Solomyak \cite{T89,KS,SoUniq}. 
studied graph-directed systems without any symmetries. They assume that there are finitely many tiles up to translation.
This approach includes the case of rotations by rational angles - rational multiples of $360^o.$ A self-similar set $A$ with pieces rotated by multiples of $60^o,$ as Figure \ref{fi2},  can be generated by a translationally finite graph-directed system of six sets without considering rotations. 

In this approach all neighbor maps are translations, and since all matrices are powers of a basic matrix $M,$ their product is commutative. 
For self-similar tilings only few examples are known outside this setting \cite{Ra94, quaquaversal, Fr08, BMT}.
For  self-similar fractals, however, we are going to present plenty of examples which are not translationally finite. 

\begin{figure}[h!] 
\begin{center}
\includegraphics[width=0.41\textwidth]{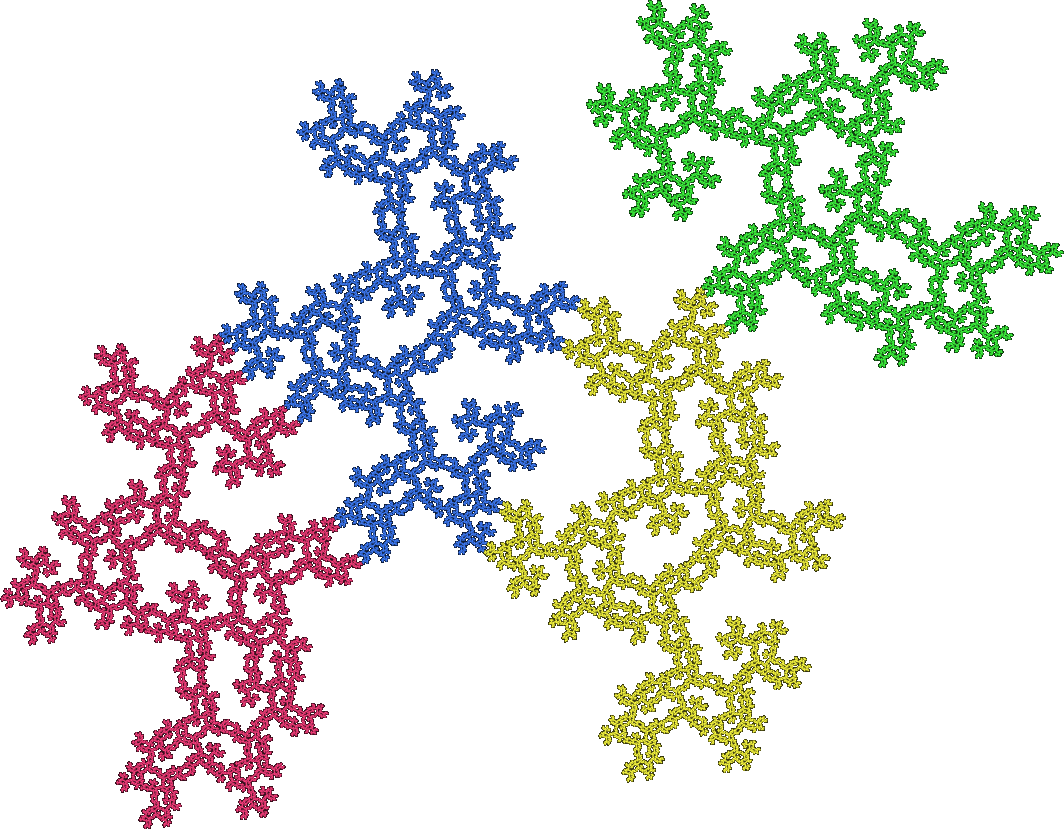} \quad
\includegraphics[width=0.54\textwidth]{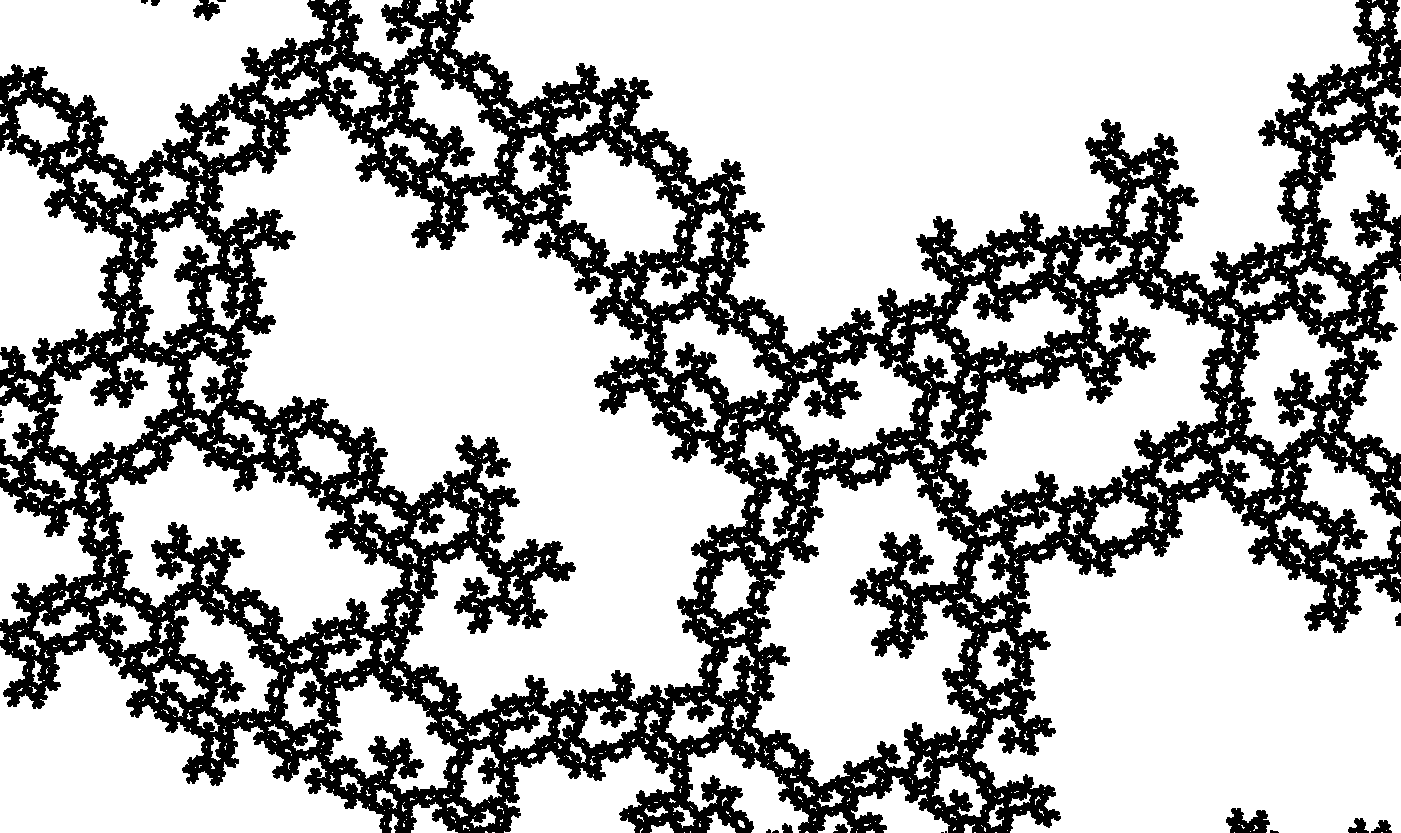}
\end{center}
\caption{A carpet generated from an irrational rotation and a reflection. It has 13 neighbor types and no characteristic directions. For details, see Section \ref{pyth}.}\label{fi7}
\end{figure}  

\subparagraph{Contents of the paper.}
In the translationally finite case as well as in self-similar sets with crystallographic data, including those which are projected from higher-dimensional lattices, all motives appear in a finite number of directions.  Here we construct fractals with a dense set of characteristic directions, in other words, with no characteristic directions at all. Figures \ref{fi6} and \ref{fi7} give a first impression of our patterns. We replace integer data by rational data, where no theorems guarantee the existence of finite type examples, except for trivial cases like Cantor sets. We performed an extensive computer search, checking some hundreds of millions of IFS, and present selected results.

After stating our basic assumptions in Section \ref{basic}, we introduce the main tool, the neighbor graph of an IFS, in Section \ref{ng}. Its topological applications are listed in Section \ref{ngapp} while Section \ref{ngalg} discusses an algebraic viewpoint and an open problem.  The search leading to Figures \ref{fi6} and \ref{fi7}, based on Pythagorean triples, is described in Section \ref{pyth}. We briefly discuss the search on the hexagonal lattice in Section \ref{hex} and introduce a classification of planar self-similar sets in Section \ref{class}.  This leads to a very convenient concept of an algebraic planar IFS which saves all matrix calculations.
In Section \ref{poly}, a family of IFS is given by a characteristic polynomial and some linear relations between expansion maps and symmetries. It turns out that the polynomial describes an algebraic number field. In the final Section \ref{quad} we discuss the search for examples in different quadratic number fields.

\subparagraph{Motivation.}
Beside the potential of isotropic fractals for modelling in science, there are various mathematical reasons for this research. One is pure curiosity: to see what is beyond crystallographic symmetries.  Another motivation is to show that our approach with symmetries is much wider than the translationally finite setting. Moreover, fractals without characteristic directions show some measure-theoretic uniformity. They have ``statistical circular symmetry'', as certain quasiperiodic tilings, and physical materials of this type show a diffraction spectrum of rings \cite{Ra94, Fr08}. 
Their projections onto lines possess the same dimension for \emph{every} direction while in general we have an exceptional set of directions with smaller dimension. Further uniformity properties were proved mainly by Shmerkin and coauthors, see \cite{FFJ15,HS12,Shmerkin15,SS18}. Here we construct concrete examples of such sets which also have nice topological properties.

\section{Basic assumptions and a simple example}\label{basic}
\subparagraph{Basic assumptions.}
We recall the basic equations \eqref{hut} and \eqref{map}:
\begin{equation} A=\bigcup_{k=1}^m f_k(A) \quad 
\mbox{ with }\quad f_k(x) = M^{-1}h_k(x) = M^{-1} (s_k x + v_k) \, .
\label{hutmap}\end{equation}
In standard coordinates, the $s_k$ should be linear isometries, given by orthogonal matrices. The map $g(x)=Mx$ should be an expanding similarity map, so that we can write equation \eqref{hutmap} as 
\begin{equation} g(A)=\bigcup_{k=1}^m h_k(A)=\bigcup_{k=1}^m s_k(A)  +v_k \, . \tag{5a}
\label{ghut}\end{equation}
Moreover, the maps must have a discrete structure to apply integer arithmetics. We assume that with respect to a common base $B=\{ b_1, b_2\}$ 
\begin{equation} M \mbox{ and the matrices of the $s_k$ contain rational numbers, and the } v_k \mbox{ are integer vectors.}
\label{intas}\end{equation}
The condition \eqref{hutmap} or \eqref{ghut} together with \eqref{intas} is our basic assumption.  
Note that $B$ need not be the standard base. It is needed to determine the combinatorial structure of the IFS. Using rational numbers with a bounded denominator, this calculation will be accurate. The transformation from base $B$ to the standard base, and the visualization of the set $A$ in standard coordinates, are done by floating point arithmetics with numerical error.  In the following sections, $B$ will be the standard base. In Sections \ref{hex}-\ref{quad} we shall consider other bases.

\subparagraph{An example with irrational rotation.}
Before we discuss the difficulties with the base $B,$ we consider examples where the standard base can be taken as $B.$ 
Let
\begin{equation}  M={\ 2\, \ 1\choose -1\ 2}\quad \mbox{ and }\quad s=\frac15 {4 -3\choose 3\quad 4}\ .
\label{ex34a}\end{equation}
The map $s$ is a rotation with irrational angle, as shown in proposition \ref{P02} below. We wonder whether an IFS with rotation $s$ between the pieces can have a nice connected attractor. A tile, as on the right of Figure \ref{fi1}, seems not possible. Thus we would like to have at least something like the Sierpi\'nski carpet. 

Let us note that a fractal tiling with irrational rotations was found in \cite{BMT} using a reflection and the fact that $M$ involves an irrational rotation. This was a rare and special example. Here we prescribe the irrational rotation directly as a symmetry.  The matrix $M$ was chosen because it has determinant 5. With two or three pieces, we have few degrees of freedom, and only a small number of IFS with OSC, most of which are well known \cite{GG,NSVW}.  With 8 or 9 pieces, we already have a huge choice of parameters $s_k$ and $v_k$ which we cannot control even with a computer (see \cite{BM19} for  discussion of a similar case). A maximal number of five pieces is better to handle. Since we find no tiles, we shall be content with fractals with $m=4$ pieces and nice topological structure.

\begin{figure}[h!] 
\begin{center}
\includegraphics[width=0.41\textwidth]{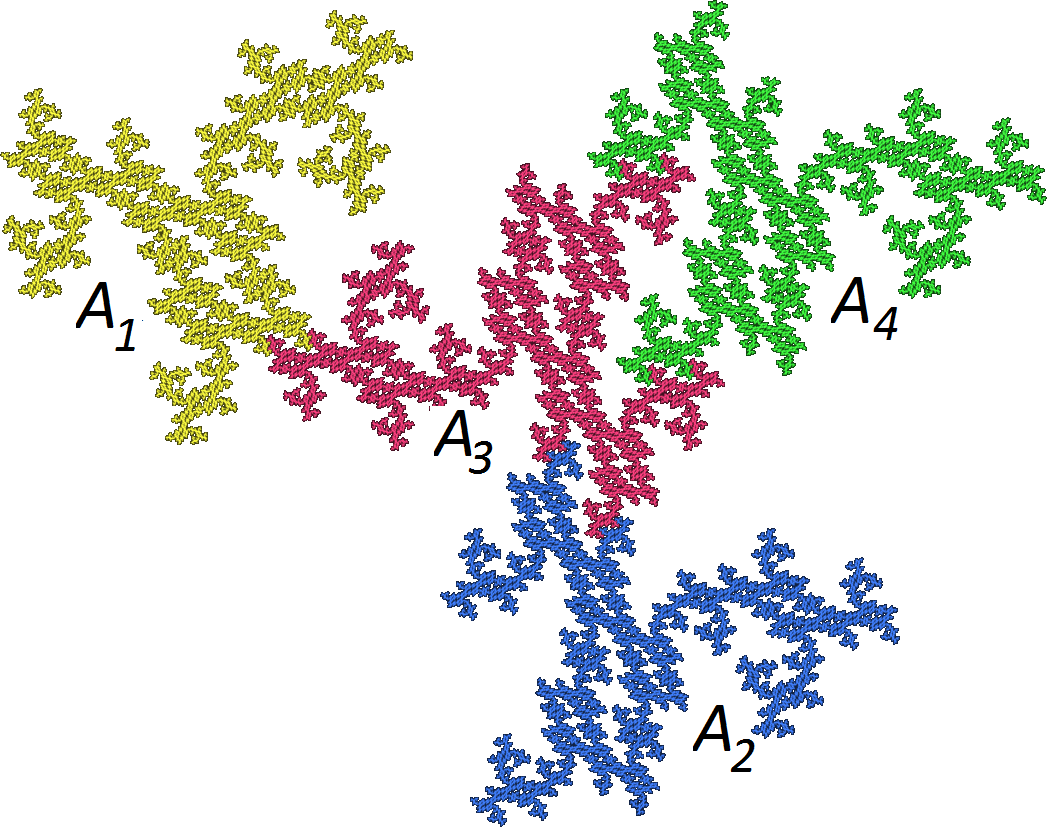}\qquad\
\includegraphics[width=0.43\textwidth]{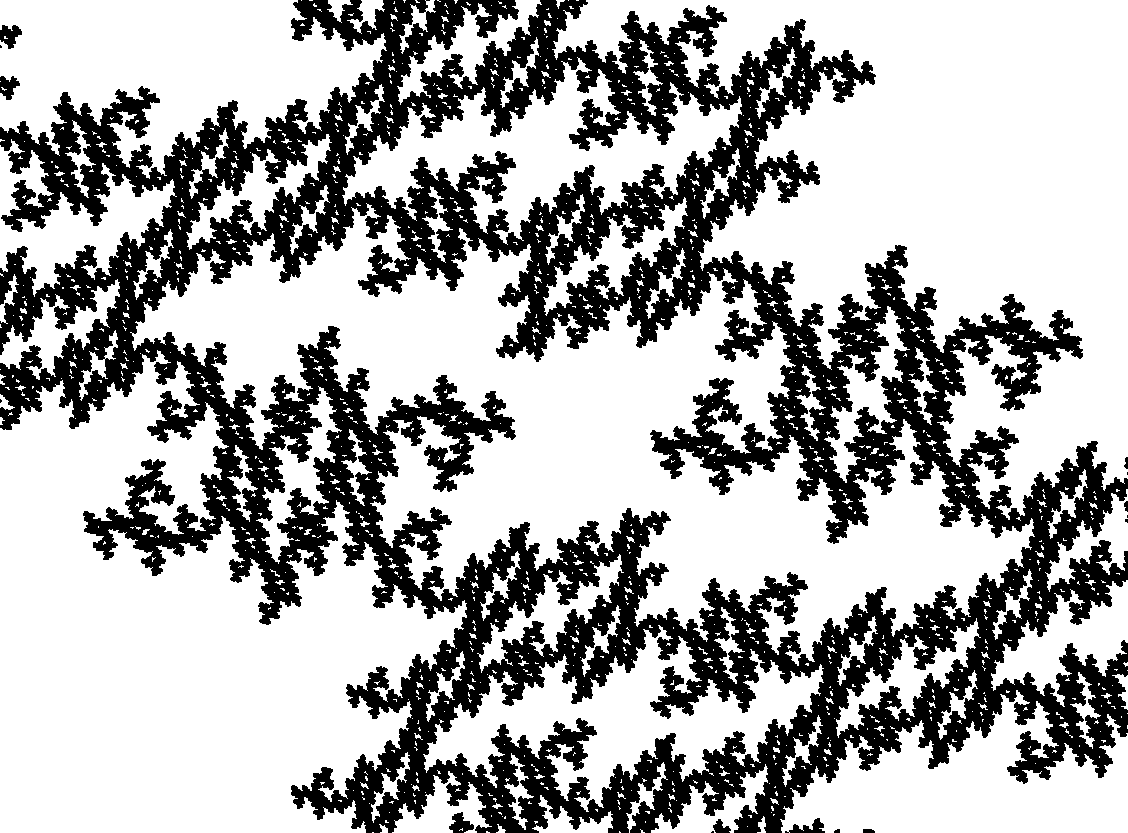}
\end{center}
\caption{A carpet with simple structure involving an irrational rotation, and a close-up.}\label{fi3}
\end{figure}  

Figure \ref{fi3} shows the simplest result. In various runs of our experiments it was always obtained early. The pieces are well-connected: they intersect in Cantor sets which all have Hausdorff dimension 0.4307, as shown below. Most importantly, there are only five neighbor types. Thus the complexity of this example is not larger as for the square or Sierpi\'nski triangle. This will be shown in the next section. First we complete the data of the IFS, listing the maps $h_k(x)=s_k(x)+v_k$ for $k=1,...,m.$

\begin{equation} \textstyle  h_1=-s(x-{0\choose 1})\ , \ h_2=-x-{0\choose 1}\ , \ h_3=x\ , \ 
h_4=-x+{1\choose 0} \ .
\label{ex34b}\end{equation}
Thus the irrational rotation acts only between the pieces $A_1$ and $A_3$ while $A_2$ and $A_4$ are parallel, and $A_3$ is rotated by $180^o$ with respect to them.

\section{The neighbor graph}\label{ng}
\subparagraph{Definition.} 
We shall now determine the combinatorial structure of our example IFS, the so-called neighbor graph. For the case of pieces of equal size, this object has been defined in various papers, including \cite{AL,Ba10,BM09,CT16,DJN12,DLN18,HLR,LZ17,RZ16,ScT,TZ20}. We introduce the concept briefly and refer to the literature for details. The IFStile package determines neighbor graphs also for IFS with different contraction ratios and for graph-directed constructions.

The neighbor graph is a directed edge-labelled graph $G=(V,E)$ which can have multiple edges and loops. The vertex set $V$ consists of all neighbor maps $h=f_w^{-1}f_v,$ where $A_u=f_u(A)$ and $A_v=f_v(A)$ are two pieces of  $A$ of the same level which intersect each other. Different pairs $(w,v)$ of words on ${1,...,m}$ can belong to the same neighbor map. 

Now if $A_w$ and $A_v$ intersect, there are intersecting subpieces $A_{wk}, A_{vj}$ which must belong to some neighbor map $h'$ which is another vertex. For each such pair of subpieces there is an edge from $h$ to $h'$ labelled by $k,j,$ indicating that
$h'=f_k^{-1}hf_j .$ Moreover, we have initial edges from $h_0=id$ to $h=f_k^{-1}f_j,$ labelled by $k,j,$  for each pair of first level pieces $A_k, A_j$ which intersect. We consider $id$ as the root vertex of $G.$ But $id$ is not a neighbor map and must not be reached by an edge when the OSC is fulfilled. So we shall not draw $id$ as a vertex. We just indicate the initial edges.

If a finite neighbor graph $G$ can be constructed for an IFS, we say that the IFS has finite type. We consider only finite type IFS. In the computer treatment, an IFS is discarded if the corresponding neighbor graph requires more than a prescribed number $N$ of vertices. In our context, we take $N=100$ or even smaller.

For our example, we first explain how a human can geometrically construct the neighbor graph. Then we describe how a computer does the job.  

\begin{figure}[h!] 
\begin{center}
\includegraphics[width=0.9\textwidth]{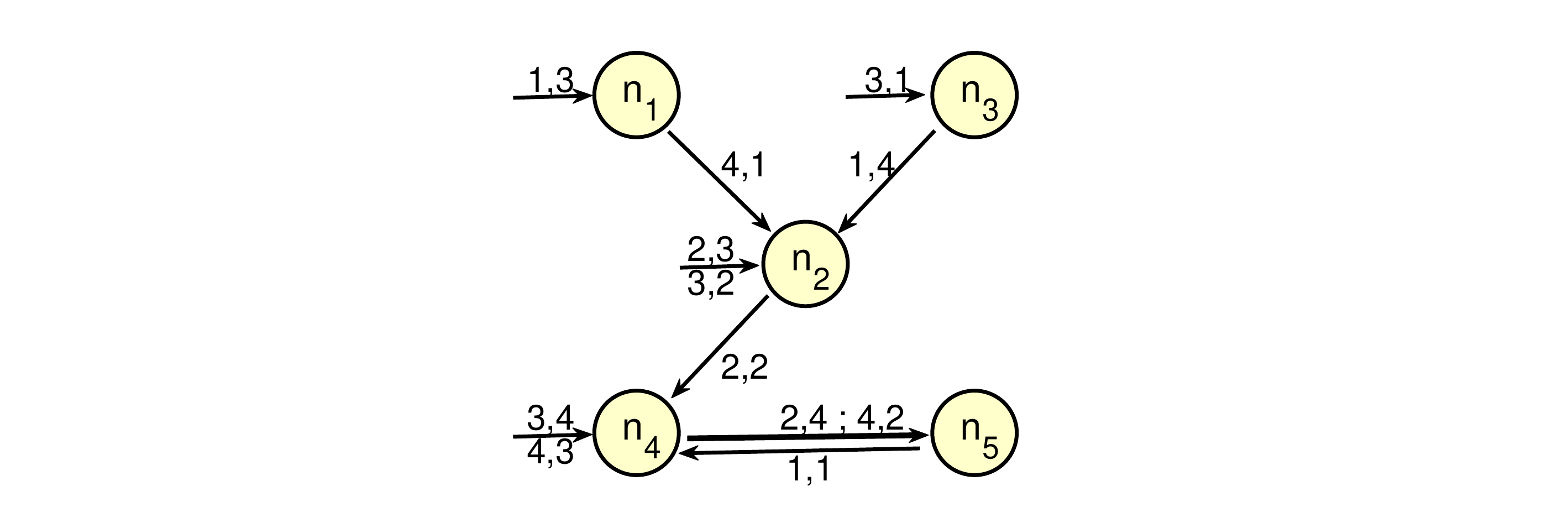}
\end{center}
\caption{The neighbor graph of Figure \ref{fi3} has a very simple structure.}\label{fi5}
\end{figure}  

\subparagraph{Intuitive construction.}
The first level intersections all involve $A_3.$ Initial edges labelled $1,3,\  3,1,\ 3,2,$ and $3,4$ go to the following vertices, respectively.
\begin{eqnarray*}\textstyle
&n_1=f_1^{-1}f_3=-s^{-1}x+{0\choose 1}, &\textstyle n_3=f_3^{-1}f_1=h_1=s(-x+{0\choose 1}),\\
&n_2=f_3^{-1}f_2=h_2=-x-{0\choose 1},  &\textstyle n_4=f_3^{-1}f_4=h_4=-x+{1\choose 0} .
\end{eqnarray*}
In this first stage, $f_k^{-1}f_j=h_k^{-1}h_j$ since $g$ cancels out, and $h_3=id$ allows to calculate directly from \eqref{ex34b}. Moreover, $n_2$ and $n_4$ are point reflections, they are self-inverse. Thus $f_2^{-1}f_3=n_2$ and $f_4^{-1}f_3=n_4.$
Instead of drawing new edges leading to the vertices $n_2$ and $n_4,$ we just write a second label $2,3$ resp. $4,3$ at the existing edges. See Figure \ref{fi5}.  

\begin{figure}[h!] 
\begin{center}
\includegraphics[width=0.42\textwidth]{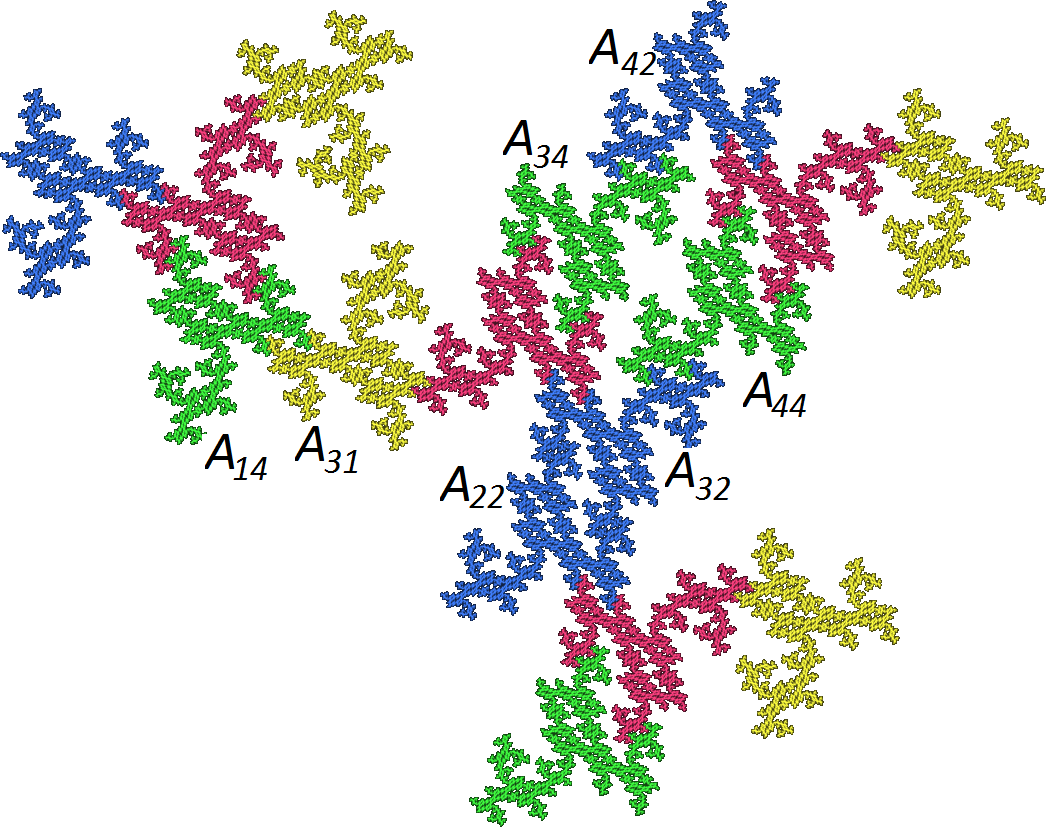}\qquad
\includegraphics[width=0.4\textwidth]{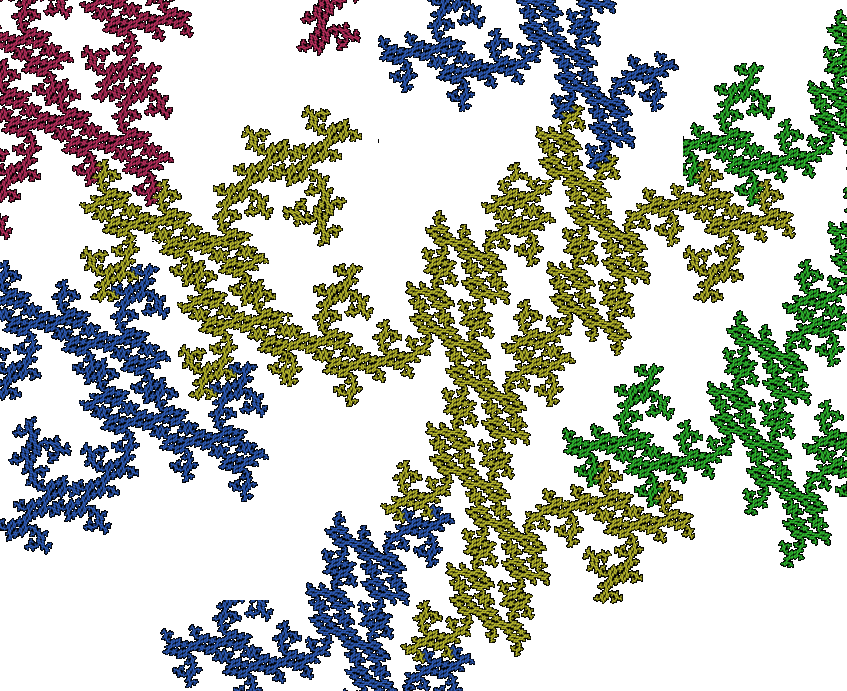}
\end{center}
\caption{Left: second level pieces Figure \ref{fi3}, colored with respect to the second index. Right: the five possible neighbors of $A$ determine the boundary sets.}\label{fi4}
\end{figure}

We have considered all first-level intersections. Now we must know how they divide into neighbor types on the second level. Figure \ref{fi4} shows our carpet, colored with respect to the second level index. We see that $A_1\cap A_3 =A_{14}\cap A_{31},$ and that the neighboring position between these subpieces is the same as between $A_2$ and $A_3.$ This can be confirmed by the equation  $f_4^{-1}n_1f_1= n_2 = f_1^{-1}n_3f_4 .$ Thus we draw an edge labelled $4,1$ from $n_1$ to $n_2,$ and an edge with label $1,4$ from $n_3$ to $n_2.$

Next, we see that $A_2\cap A_3 =A_{22}\cap A_{32},$ and the neighbor type of these subpieces is just $n_4.$ So we draw an edge with label $2,2$ from $n_2$ to $n_4.$
Moreover, the neighbor type $n_4$ divides into two subtypes:
\[ A_3\cap A_4 =(A_{32}\cap A_{44}) \cup (A_{34}\cap A_{42}) \ . \]
The corresponding neighbor map   
\[\textstyle  n_5=f_2^{-1}n_4f_4=-x-{1\choose 0} \] 
agrees with $f_4^{-1}n_4f_2.$ So we draw an edge with two labels $2,4$ and $4,2$ from $n_4$ to $n_5.$  The subdivision of neighbor types  $n_1$ to $n_4$ is completed. We still have to study the subdivision of the new type $n_5$ by considering the third level pieces of $A,$ Fortunately, it can be seen that $n_5$ describes just the intersection of the subpieces with index 1 on both sides, and represents type $n_4$ on the next level. Thus only one more edge from $n_5$ to $n_4$ with label $1,1$ has to be drawn, resulting in Figure \ref{fi5} which is discussed in the next section.

\subparagraph{The computer algorithm.}
A computer easily generates lots of neighbor maps  $h_{wv}=f_w^{-1}f_v$ by repeatedly applying the recursive formula $h'=f_k^{-1}hf_j$ with $k,j\in\{ 1,...,m\}.$  The problem is to decide which of these actually fulfil  $f_w(A)\cap f_v(A)\not=\emptyset .$  

\begin{Proposition}\label{P0}
For closed sets $A,B$ let $d(A,B)=\inf \{ |x-y|\, | x\in A, y\in B\} .$ Suppose that an isometry $h$ fulfils $d(h(A),A)\ge\varepsilon >0.$ Then  $d(h'(A),A)\ge\varepsilon  /r$ for each $h'=f_k^{-1}hf_j$ with $k,j\in\{ 1,...,m\}.$ 
\end{Proposition}

\emph{Proof. } We have $d(h(f_j(A)),f_k(A))\ge d(h(A),A)\ge\varepsilon .$ Applying the similarity map $f_k^{-1}$ with factor $1/r$ to both sets, we conclude $d(
f_k^{-1}hf_j(A),A) \ge\varepsilon /r.$
\hfill $\Box$ \vspace{2ex}

The proposition says that when  $h(A)$ does not intersect $A$ for some generated map $h,$ then $h'(A)$ will not intersect $A$ for all its successor maps $h'$ in forthcoming levels. Moreover, the size of the translation of $h'$ will grow exponentially with the level. So even for small $\varepsilon$ it will become obvious  after few recursion steps that we do not have a neighbor map.

The algorithm is now clear. We stop the recursive calculation as soon as we recognize that $h'$ is not a neighbor map. Afterwards, we repeatedly remove from our recursive tree of maps all vertices $h$ without successors.  Then we are left with all maps which lead to a cycle in the neighbor graph, and this is exactly the graph $G.$  On the other hand, when the recursive calculation produces too many (say $10^5$) isometries for which we cannot decide whether they are neighbor maps, we give up and say that the IFS seems not to be of small finite type.

Now we provide a simple general criterion to decide when an isometry is not a neighbor map. The IFStile package uses more complicated IFS-specific estimates which reduce the effort of recursive calculations.

\begin{Proposition}\label{P01}
Let $\tilde{x}$ denote the center of gravity of $A,$  determined from the IFS by the equation $\tilde{x}=\frac1m \sum_{k=1}^m f_k(\tilde{x}) \ .$ Moreover, let 
$\delta =\max_{k=1}^m |f_k(\tilde{x})-\tilde{x}| \ .$ \\
Then $A$ is contained in the ball $B$ with radius $\delta /(1-r)$ around 
$\tilde{x}.$\\
When an isometry $h$ fulfils  $|h(\tilde{x})-\tilde{x}|> 2\delta /(1-r)$ then $h(A)$ and $A$ are disjoint.
\end{Proposition}

\emph{Proof. } Note that $\tilde{x}$ and $\delta$ can be easily calculated from the IFS data. For each $j,k$ the distance of $x=f_jf_k(\tilde{x})$ to $f_j(\tilde{x})$ is at most $r\delta .$ So $|x-\tilde{x}|\le (1+r)\delta .$ For a word $w=k_1...k_n$ and $x= f_w(\tilde{x})$ we get $|x-\tilde{x}|\le (1+r+...+r^{n-1})\delta .$ Since these points, with arbitrary $n$ and $w,$ are dense in $A,$ we obtain the estimate
$|y-\tilde{x}|\le \delta /(1-r)$  for all $y\in A.$ This proves the first assertion. If $|h(\tilde{x})-\tilde{x}|$ is greater than twice the radius of a ball $B$ with centre $\tilde{x},$ then $B$ and $h(B)$ must be disjoint. 
\hfill $\Box$ \vspace{2ex}

\section{What can we conclude from the neighbor graph?} \label{ngapp} 
\subparagraph{Topology-generating automaton.}
The neighbor graph is an automaton which generates the topology of $A.$ \emph{All} topological properties of $A$ are encoded in the neighbor graph. In the argument above, we pretended that we can see that the pieces $A_1$ and $A_3$ intersect. This was not true. Actually, pictures did repeatedly lead us to wrong conclusions. Only the calculation of the neighbor graph can verify that $A_1$ and $A_3$ have common points.

The neighbor graph is empty if the pieces $A_k=f_k(A), \ k=1,...,m$ are pairwise disjoint. Otherwise, the neighbor graph tells us which pieces intersect. Consider any infinite path of edges, starting with an initial edge, that is, starting from the root  $h_0=id.$ Paths are directed by definition, and if they are infinite, they must contain directed cycles.  If the path has the labels $k_1,j_1; k_2,j_2; ...$ then the two sequences 
$k_1k_2...$ and $j_1j_2...$ describe the same point, which is determined by the decreasing sequence of compact sets $A_{k_1}\cap A_{j_1},  A_{k_1k_2}\cap A_{j_1j_2}$ and so on (cf. \cite{Bar,BP,Fal}). In our example, infinite paths are all given by the cycle between $n_4$ and $n_5$ which we have to run through infinitely many times. Actually, it consists of two cycles since from $n_4$ to $n_5$ we have two edges with labels $2,4$ and $4,2,$ respectively. Thus the addresses $3u_11u_21u_3...$ and $4v_11v_21v_3...$ describe the same point whenever either $(u_i,v_i)=(2,4)$ or  $(u_i,v_i)=(4,2)$ for each $i=1,2,... .$ This means that $A_3$ and $A_4$ have a whole Cantor set in common. And the same holds for $A_3$ and $A_2,$ and $A_3$ and $A_1,$ if we consider paths starting in $n_2$ and $n_1.$ On the other hand, $A_1$ and $A_4$ have no common points since $1,4$ and $4,1$ do not appear as labelling of an initial edge.

\subparagraph{Connectedness properties.}
Recall that the neighbor graph $G$ is constructed so that there is an outgoing edge from each vertex, and hence an infinite directed path starting in each vertex.  Thus there is an initial edge with label $k,j$ if and only if $A_k\cap A_j\not=\emptyset .$
So from $G$ we can define the \emph{connectedness graph} $G_c$ of $A,$ with vertex set $\{ 1,...,m\}$ and undirected edges between $k,j$ whenever $A_k$ intersects $A_j.$  For our example, the connectedness graph is a 3-star, or letter Y, with central vertex 3. The following is known. The first assertion is a classical theorem of Hata,
proved by inductively constructing chains of intersecting pieces.  

\begin{Proposition}\label{PC}
\begin{enumerate}
\item[(i)] The attractor $A$ is connected if and only if $G_c$ is connected.
\item[(ii)]  Suppose $A$ is connected. Then $A$ contains a closed Jordan curve if and if and only if either $G_c$ is connected, or two pieces  $A_j, A_k$ intersect in more than one point.
\item[(iii)]  Let $h=f_k^{-1}f_j$ be the vertex of the neighbor graph $G$ which corresponds to $D=A_k\cap A_j\not=\emptyset .$ Let $\cal C$ denote the set of all directed cycles which can be reached by a path from vertex $h.$ 
The intersection $D$ is a singleton if and only if $\cal C$ consists of one element $C_1$ and there is only one path from $h$ to $C_1.$ The intersection is finite if there is no path between $C_1$ and $C_2$ for any $C_1, C_2\in \cal C .$  The set $D$ is uncountable if there exist $C_1, C_2\in \cal C$ and paths from $C_1$ to $C_2$ and back. Otherwise $D$ is countably infinite. 
\end{enumerate}
\end{Proposition}

The proof is simple but there are various details, cf. \cite{BM09}.
The graph in Figure \ref{fi5}, for instance, contains two cycles from $n_4$ through $n_5$ to $n_4.$ The two cycles can be reached from each other since they start in the same point - the connecting path is empty.  Similarly, care has to be taken when points belong to three or more pieces. 

According to (iii), we have different levels of connectedness of $A$ which we can read from the neighbor graph: ordinary connectedness in the sense of Hata,
connectedness by single-point intersections which means that $A$ is a dendrite, connectedness by finite intersections (called p.c.f. self-similar sets in \cite{Kig}), connectedness by infinite, at most countably infinite or uncountable intersections.  All such properties will be determined in IFStile, and can be used to  characterize and order the exmples.

Some related properties still require research, for instance the study of connectedness components of $A$ when $A$ is not connected. For distinguishing Cantor set intersections $A_k\cap A_j$ from intervals, one needs to consider neighbor graphs for intersections of three or more pieces \cite{Ba10,CT16,TZ20}. This is important for tilings, especially in three dimensions. In the example above, all three-piece intersections are empty.

\subparagraph{Cutpoints and first-level intersections.}
Consider the fixed point $y=0$ of $f_3$ in the above example. It does not belong to two pieces since the address $\overline{3}=333...$ - even the word 33 - is not labelling a path in the neighbor graph. Nevertheless, $y$ is an important point for the topology of $A.$ It is a \emph{cutpoint} of degree 3: $A\setminus\{ y\}$ consists of three connected components.  This can be easily seen: removal of $A_3$ from $A$ results in three different pieces, and removal of $A_{33}$ from $A$  results in three larger pieces, and so on.  Thus $A$ involves a rather simple tree structure, despite containing Cantor set intersections and closed Jordan curves.
We call such a space a web.

We do not go further into detail. We just want to give an impression of the potential and universal role of the neighbor graph. For the computer search, certain simple parameters are useful even if they do not define topological properties of $A.$ The number FLI of first-level intersections is just the number of initial edge labels $k,j$ with $k<j.$ This is an invariant of the IFS (cf. \cite{Sierrel}). In our example ${\rm FLI}=3.$  When we look for connected examples, and do not find them - which was the case in several experiments related to this paper - then FLI will give us an idea on how far we are from connectedness.

\subparagraph{Hausdorff dimension of the boundary.}
Apart from connectedness properties, let us remind the fact that the neighbor graph verifies that our IFS is of finite type and fulfils the OSC. In particular, the Hausdorff dimension of $A$ equals \cite{Bar,BP,Fal}
\begin{equation} \alpha =\frac{\log m}{-\log r}= \frac{2\log m}{\log \det M}\ .
\label{simdim}\end{equation}
Moreover, $A$ is a very nice type of self-similar set, completely described by an automaton. It is a computable self-similar set.

Now we show that the boundary dimensions of $A$ can also be read from the neighbor graph. This was found for the twindragon by Gilbert \cite{Gi86}, for the L\' evy dragon by \cite{DKV,SW} and for self-affine tiles by \cite{AL,DJN12,LL13} and others. To each neighbor type $n_k$ there corresponds a boundary set $B_k=A\cap n_k(A).$ The first parts of the labels of the edges are used to determine a set of equations for the boundary sets:
\[ B_2=f_4(B_1)=f_1(B_3), B_4=f_2(B_2)=f_1(B_5), B_5=f_2(B_4)\cup f_4(B_4)\ .\]
The boundary sets form a graph-directed system, which inherits the OSC from $A.$  Therefore it is possible to determine the Hausdorff dimensions of all boundary sets. In our case they all have the same dimension, and we see that $B_4$ is even a self-similar set:
\[  B_4=  f_1f_2(B_4) \cup f_1f_4(B_4) \ .\]
The 4 mappings $f_k$ have the factor $r=1/\sqrt{5}.$ So $A$ has dimension
\[ \alpha =\frac{\log 4}{\log\sqrt{5}}=\frac{4\log 2}{\log 5}\approx 1.7227\ ,\]
and all the boundary sets have dimension   
\[ \beta =\frac{\log 2}{\log 5}=\frac{\alpha}{4}\approx 0.4307\ .\]
Note that with $f_1,$ the irrational rotation is involved in all boundary sets.

All the mentioned properties and parameters are automatically determined by IFStile, and can be used to identify non-isomorphic examples, or to select specimen with prescribed topological properties from a computer-generated collection, cf. \cite[Section VI]{Sierrel}.

\section{Algebraic aspects of the neighbor graph}\label{ngalg}
\subparagraph{An open question.} The neighbor graph completely describes the topology of $A,$ but not the Hausdorff dimension of $A.$ It is well-known that the dimension of a Koch curve can be changed without changing the topology, just by varying the angle of the maps \cite{Bar}. So the question is which properties   should be added to the neighbor graph in order to completely characterize the IFS, up to isomorphy. Isomorphy involves change of the coordinate system and permutation of the maps, see \cite{Sierrel}. The question seems difficult. We shall discuss it for our simple example.

\subparagraph{Generating relations.}
The neighbor maps generate a group of isometries. Whenever we extend the self-similar construction of a connected attractor $A$ to the outside, by forming supertiles $f_{k_1}^{-1}(A),\  f_{k_2}^{-1}f_{k_1}^{-1}(A),...,$ any isometry between two `tiles'  of such a pattern will belong to that group. In algebra, groups are often defined by a system of generators and generating relations. It turns out that the neighbor graph provides such relations between the neighbor maps. In this way it gives algebraic information on the IFS beyond the topology of $A.$ 
Usually the equations are highly nonlinear. For our simple example we can separate and solve them, however.

\subparagraph{Finding IFS data from a neighbor graph.}
Let us assume that we have an IFS given by $g,h_1,...,h_4$ which produces the neighbor graph of Figure \ref{fi5}. For simplicity we assume $h_3=id$ which can  be arranged by passing to the isomorphic IFS $\tilde{g}=h_3^{-1}g, \ \tilde{h_j}=h_3^{-1}h_j.$  Then the initial edges 3,1, \ 3,2 and 3,4 imply that $h_1=n_1, \ h_2=n_2,$ and $h_4=n_4.$ 

If two initial edges with label $k,j$ and reversed label $j,k$ lead to a vertex $h,$ then the map $h$ is self-inverse. This also holds for paths with reversed labels starting from the root. Thus
$n_2, n_4,$ and $n_5$ are self-inverse. Now when we restrict ourselves to IFS without reflections, the only self-inverse isometries in the plane are point reflections. Thus we can conclude that $h_2(z)=w_2-x, \ h_4(x)=w_4-x,$ and $n_5(x)=w_5-x,$ for some vectors $w_k.$

In order to simplify equations and avoid confusion with the IFS data above, we now work in the complex plane, writing $z$ instead of $x.$ Our aim is to determine unknown complex numbers $\lambda , t, $ and $v$ from the neighbor graph. They correspond to $M, -s,$ and $s{0\choose 1}$ in the example. The $w_k$ are now also considered unknown complex numbers, but we use the same letter as for vectors above, and keep the notation $g,h_k,f_k$ for the unknown maps. We have proved the first part of (i) in the following proposition.

\begin{Proposition}\label{PA}
Suppose an IFS is defined in the complex plane by $g(A)=\bigcup_{k=1}^4 h_k(A)$ with $g(z)=\lambda z,\ h_1(z)=tz+v,\ h_3(z)=z,$ with $\lambda,t,v \in\CC$ and $|t|=1.$ If the neighbor graph of this IFS is the graph in Figure \ref{fi5} then 
\begin{enumerate}
\item[(i)] $h_2(z)=w_2-z$ and $h_4(z)=w_4-z$ where $w_2,w_4\in\CC$ fulfil $w_4=(2-\lambda)w_2.$
\item[(ii)]  $v=tw_2.$
\item[(iii)]  $  t(4-\lambda) =\lambda (\lambda^2 -3\lambda +3)\ .$
\end{enumerate}
\end{Proposition}

\emph{Proof. } For (i) we derive the relation between $w_2$ and $w_4$ from the neighbor graph. Note that $f_k=\lambda^{-1}h_k.$ The edge from $n_2$ to $n_4$ generates the relation $f_2^{-1}n_2f_2=n_4.$ We multiply by $f_2$ from the left and use $h_2=n_2,\ h_4=n_4.$ Thus 
\[ h_2\lambda^{-1}h_2=\lambda^{-1}h_2h_4,\quad\mbox{ or }\
w_2-(\lambda^{-1}(w_2-z)=\lambda^{-1}(w_2-(w_4-z))\, .
\] 
Ordering terms, we obtain (i).\\
The edge from $n_1=h_1$ to $n_2=h_2$ implies $h_1f_4=f_1h_2,$ or
$t(\lambda^{-1}(w_4-z)+v)=\lambda^{-1}(t(w_2-z)+v).$ Thus $tw_4+\lambda v= tw_2+v.$ Together with the relation in (i), this gives (ii).\\
To verify (iii), consider the edge from $n_4$ to $n_5.$ It generates the relation $f_2^{-1}n_4f_4=n_5,$ or $h_4\lambda^{-1}h_4=\lambda^{-1}h_2n_5.$ Calculation gives
\[ w_5=w_4-\lambda w_4+w_2= (\lambda^2 -3\lambda +3)w_2\ .\]
The edge from $n_5$ to $n_4$ says
$n_5f_1=f_1n_4,$ or $w_5-\lambda^{-1}(tz+v)= \lambda^{-1}(t(w_4-z)+v).$
Thus
\[ \lambda w_5= tw_4 +2v= t(w_4+2w_2)=t(4-\lambda) w_2 \ .\]
Comparing the two expressions for $w_5$ we get (iii).  
\hfill $\Box$ \vspace{2ex}

The parameter $w_2$ was not determined. Since the assumption $h_3=id$ did not involve coordinate axes, we are still free to choose the coordinate system.
We can set $w_2=1,$ or $w_2=-i$ if we like to obtain equation \eqref{ex34b} exactly.  Thus we have almost determined the whole IFS from the neighbor graph. We got a rational function $t=r(\lambda)=\lambda (\lambda^2 -3\lambda +3)/ (4-\lambda).$ For $\lambda_0=2-i,$ the value $t_0=r(2-i)=(2-i)(-i)/(2+i)=-(4+3i)/5$ represents the symmetry $s$ in \eqref{ex34a}.

Since $t$ varies on the unit circle, it seems better to consider $\lambda$ as a function of $t,$ that is, of the rotation angle of $h_1.$ The derivative of $r$ is $r'(\lambda)=(3(\lambda -1)^2+r(\lambda))/(4-\lambda ).$ Since $r'(2-i)=(-6i+t_0)/(2+i)\not= 0,$ the inverse function theorem applies: there is a holomorphic function $\lambda =r^{-1}(t)$ defined in a neighborhood of $t_0.$ Moreover, in a sufficiently small neighborhood of $t_0,$
 the maps $f_k^{-1}f_j$ which correspond to non-intersecting pieces at $(t_0,\lambda_0)$ will still correspond to non-intersecting pieces at $(t,\lambda=r^{-1}(t))$ since only a finite number of such pairs need to be considered. And the above calculations can be reversed so that the parameters $(t,\lambda=r^{-1}(t))$ provide IFS with the neighbor graph of Figure \ref{fi5}.  This proves

\begin{Proposition}\label{PB}
In a small neighborhood of our example IFS parameters $\lambda_0=2-i, \ t_0=-(4+3i)/5$ together with $w_2=1, v=t, w_4=2-\lambda,$ there is a one-dimensional parametric family of IFS which all have the same neighbor graph as the example. In this neighborhood and family, each IFS is uniquely determined by the neighbor graph and $|\lambda|$ which respresents Hausdorff dimension of the attractor.
\end{Proposition}

The last assertion was checked only by a Matlab calculation which shows that $|\lambda|$ is decreasing aolmost linearly on the curve $\lambda =r^{-1}(t)$ for $30^o\le {\rm angle} t\le 45^o$ while at $t_0$ we have the angle $\arctan\frac34\approx 36.9^o.$ 

Since (iii) is a cubic equation, there are two other solutions $\lambda_1, \lambda_2$ which fulfil $r(\lambda_j)=t_0.$  One of them has modulus smaller than one, and for the other one IFStile could not calculate a neighbor graph, so it certainly will not generate Figure \ref{fi5}. Thus it seems even globally true that for this special neighbor graph, an associated IFS is uniquely determined by its Hausdorff dimension.

\section{Examples from Pythagorean triples}\label{pyth}
\subparagraph{Pythagorean triples and irrational rotation.} 
A rotation of the plane by an angle $\alpha$ is called rational if $\alpha=\frac{k}{n}\cdot\pi$ for integers $k,n,$ and irrational otherwise. 

\begin{Proposition}\label{P02}
A rotation $s$ has a matrix of the form ${u -v\choose v \ u}$ with $u^2+v^2=1.$  If $u,v$ are rational then the rotation angle is irrational or a multiple of $90^o.$ 
\end{Proposition}

\emph{Proof. } If there are integers $k,n$ with $\alpha=\frac{k}{n}\cdot\pi$ then
$z=u+iv$ is an $n$-th root of unity in $\CC .$ That is, $z$ is a root of the $n$-th  cyclotomic polynomial which is irreducible. Thus except for $n=1,2,$ and 4, $z$ is not in the rational field $\QQ(i),$ and $u$ or $v$ must be irrational.
\hfill $\Box$ \vspace{2ex}

The rational points $(u=a/c, v=b/c)$ on the unit sphere correspond to the Pythagorean triples $(a,b,c)$ of integers with $a^2+b^2=c^2.$ They can be generated by a well-known formula of Euclid, see Section \ref{quad}.
The rotation $s$ in \eqref{ex34a} came from the basic triple $(3,4,5).$ The two next triples are $(5,12,13)$ and $(8,15,17).$

Our question here is whether the corresponding irrational rotations, combined with appropriate expansions $g(x)=Mx,$ do generate self-similar sets with strong connectedness properties. 

\begin{figure}[h!] 
\begin{center}
\includegraphics[width=0.38\textwidth]{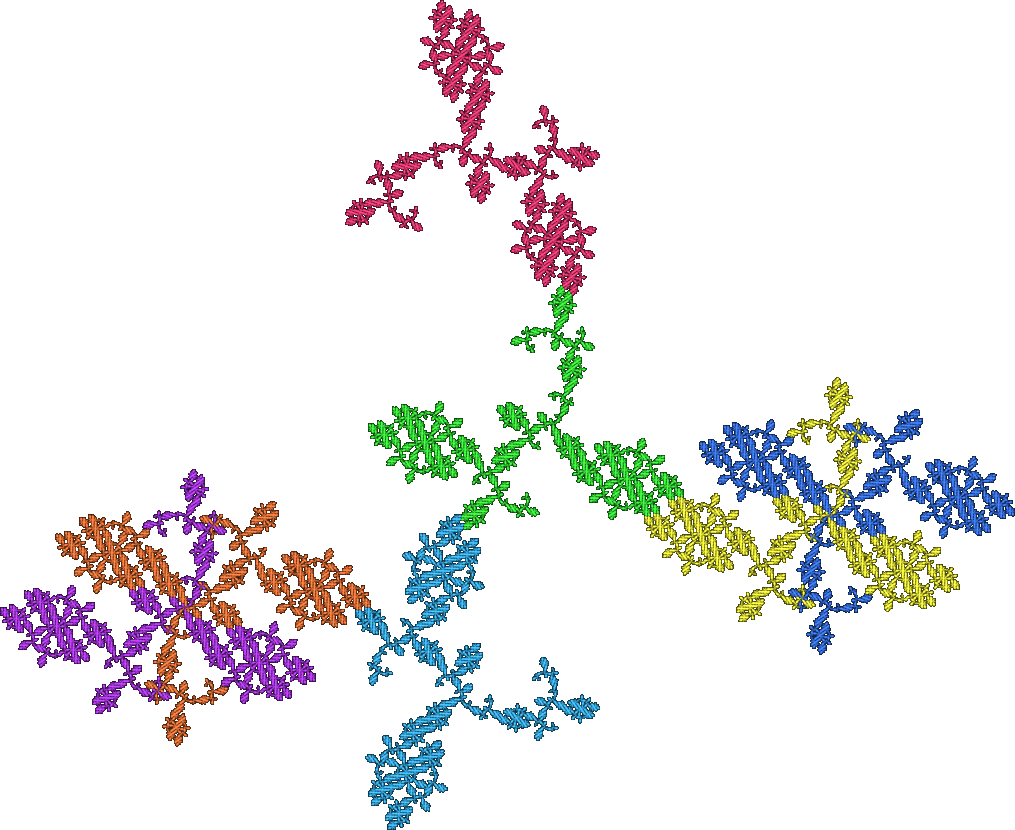} \qquad
\includegraphics[width=0.5\textwidth]{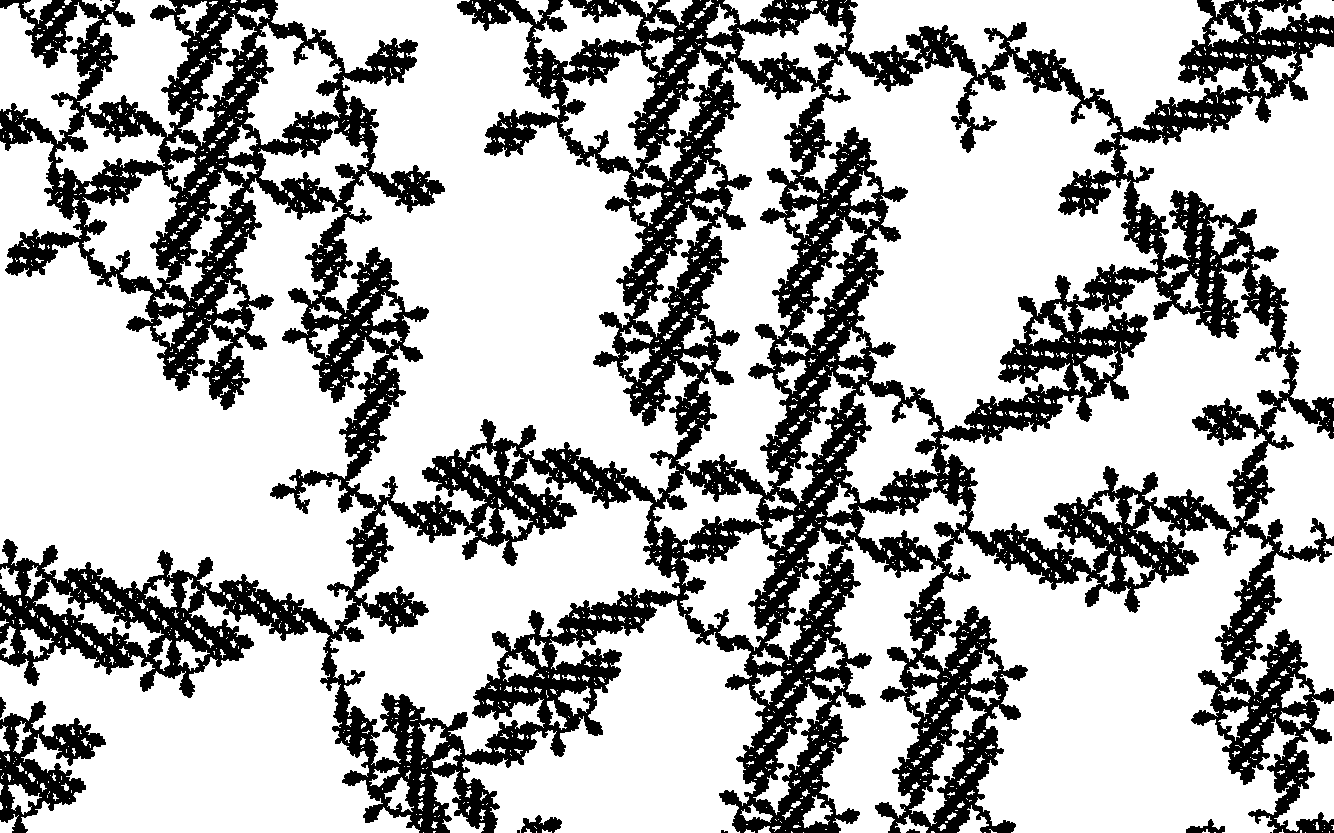}
\end{center}
\caption{A pattern of dimension 1.69 generated from $(3,4,5)$ together with $90^o$ rotations. Such patterns are fairly easy to generate. This one has only 9 neighbor types.}\label{fi8}
\end{figure}  

\subparagraph{The challenge.} 
We briefly explain the goal of our experiments. For crystallographic data, each integer expansion matrix $M$ leads to IFS with OSC and $m=\det M.$ This is proved by taking so-called complete residue systems as $h_1,...,h_m$ \cite{Ba5}, and it is well-known that the attractors are tiles. By dropping one or more of the mappings we can then easily create carpets, as in Figure \ref{fi1}.
For irrational data, however, only very few tiles are known \cite{Ra94,quaquaversal,Fr08,BMT}, and we found no further tiles in our experiments. 

It is not difficult 
to find totally disconnected, Cantor-type sets $A$ with all kinds of irrational IFS.  But whenever two pieces have a point in common, this creates an equation for the mappings of the IFS. In general the equation is given as a limit when the level of pieces tends to infinity. However, if we assume finite type, edges of the neighbor graph define  equations for finite compositions of the $f_k,$ as we have seen above. In case of complex linear functions $f_k(z)=t_kz+v_k$ we directly get polynomial equations in $z.$

When we require the pieces to have only one point with eventually periodic addresses in common (so-called p.c.f. fractals, like the Sierpi\'nski triangle \cite{Kig,Str}), we can still establish such equations by hand and try to solve them. 
However, we do not know any mathematical arguments which would lead to the construction of the example above, even though there is only one double cycle in the neighbor graph. It seems a mystery that two pieces differing by an irrational rotation intersect in a similar Cantor set as pieces which just differ by a point reflection.

\subparagraph{Results of the computer experiments.} 
While a computer search with rational rotations gives thousands of examples within a minute, many of them with high complexity \cite{Sierrel}, the interactive search for examples with irrational IFS is more difficult. First we get only Cantor sets without any intersections. Then few examples may have pieces with intersections, creating a non-trivial neighbor graph, but $A$ will still remain a Cantor set. Then by modifying the best datasets, we get some fractals which are connected, or have at least connected subsets. Finally, deleting all bad examples and modifying only the best ones, we get carpets with uncountable intersections in most of our experiments. As a rule, they have small complexity, between 10 and 40 neighbor types. The number of such carpets varied between 1 and 100, for different $g$ and $s.$

Two results were given in Figures \ref{fi6} and \ref{fi7} in Section \ref{intro}.
They are from the same family as Figure \ref{fi3}, with $g$ and $s$ defined in \eqref{ex34a} and $m=4.$ Hence they have the same Hausdorff dimension $\alpha \approx 1.72.$
We shall not further specify the IFS data since all .png files produced by IFStile contain their data. They can be opened in IFStile, allowing visual study of details.  The IFS data, including code of the neighbor graph, can be found under 'View-Console'. A list containing all our figures will be available on the web page \cite{M}.

Most of our examples with irrational rotations contained cutpoints like Figure \ref{fi3}. We thought that such cutpoints exist for mathematical reasons. Figure \ref{fi6}, with 16 neighbor types and boundary dimension 0.24, seems to reject this conjecture. It has a single point intersection of two pieces. However, this is not a global cutpoint. It seems that there are only local cutpoints which separate a small neighborhood.

If a reflection is allowed as additional symmetry, we get lots of carpets like Figure \ref{fi7}. They have no cutpoints, neither global nor local ones. The influence of the reflection is visible in the local structure. In contrast to Sierpi\'nski's triangle and carpet, our examples do not contain line segments. It would be interesting to know how much fractal analysis on these spaces differs from the classical carpet. 

\begin{figure}[h!] 
\begin{center}
\includegraphics[width=0.35\textwidth]{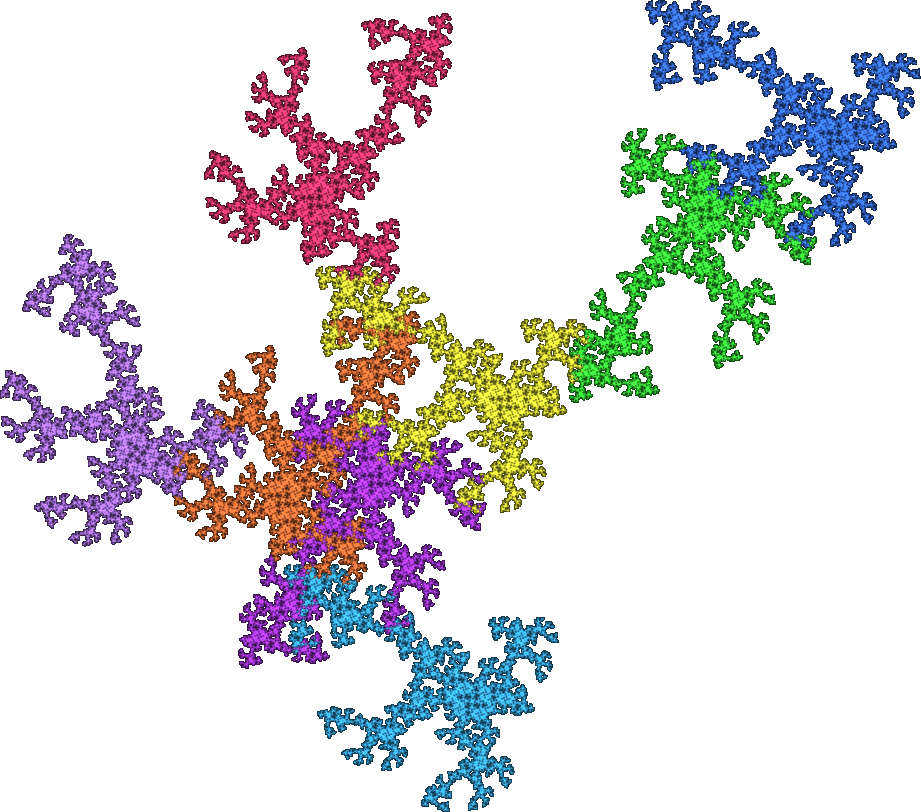} \qquad
\includegraphics[width=0.49\textwidth]{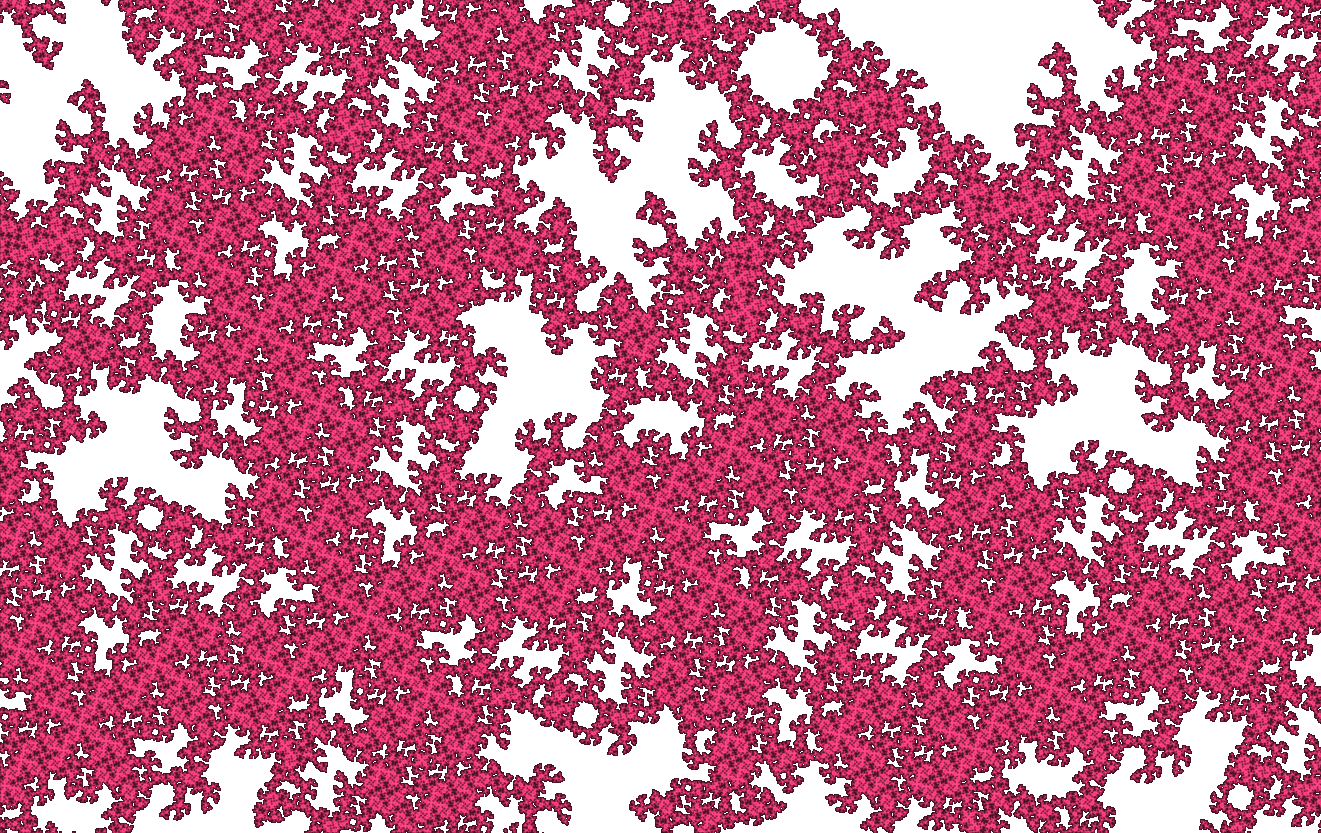}
\end{center}
\caption{A carpet generated from $(3,4,5)$ together with $90^o$ rotations. The expansion has determinant 10, and there are $m=8$ pieces. The Hausdorff dimension is 1.81, the intersections of pieces have dimension 1.25, and there are 13 neighbor types. }\label{fi9}
\end{figure}  

Any additional symmetry will increase the number of patterns.
Instead of a reflection, a rotation by $90^o$ can be taken. 
A similar effect arises when we take two Pythagorean rotations with angles adding to $90^o.$ We tried this for the expansion $g=5s-3$ with matrix
${1 -3 \choose 3\ 1 }$ which has determinant 10. A triangular tile with $m=10$ and  a dense set of characteristic directions is known \cite{Encyc}. It uses a reflection, as in \cite{Ra94,BMT}. For the Pythagorean rotation together with $90^o$ rotations we found plenty of patterns with $m=7$ like Figure \ref{fi8}, with $m=8$ like Figure \ref{fi9}, and various carpets with $m=9$ and Hausdorff dimension 1.91.  The influence of the $90^o$ rotations often dominated the effect of the irrational rotation.  We found a few similar figures with other Pythagorean triples.  Most examples had rather small complexity, in contrast to rational examples like Figure \ref{fi2}. 

\section{An example from the hexagonal lattice}\label{hex}
If we want to include $60^o$ rotations, we have to distinguish the standard base  and the base $B=\{ b_1,b_2\}$ mentioned in Section \ref{basic} for which the matrices of $g$ and $s_k$ have rational entries.
A counterclockwise rotation $s$ by $60^o$ has matrix $0\ -1\choose 1\ \ 1$ with respect to the vectors $b_1={1\choose 0}, \ b_2={1/2\choose \sqrt{3}/2} .$ When we work with such rotations, as in Figure \ref{fi2}, we determine the neighbor graph of the IFS with respect to the base $B,$ using accurate integer calculation. For visualization on a numerical level we transform back to standard coordinates. 

As an example, we take $g=2s+1$ as expansion and $t=(3s+5)/7$ as irrational rotation. Here 1 denotes the identity map. Thus $g(x)=2s(x)+x$ is the similarity map transforming $b_1$ into $2b_2+b_1={2\choose\sqrt{3}},$ with determinant 7. By expressing $g$ and $t$ as rational linear combinations of $s$ and $1,$ we guarantee that their matrices with respect to base $B$ are rational, and that the mappings commute. Moreover, we save a lot of matrix calculations. 
 
\begin{figure}[h!] 
\begin{center}
\includegraphics[width=0.35\textwidth]{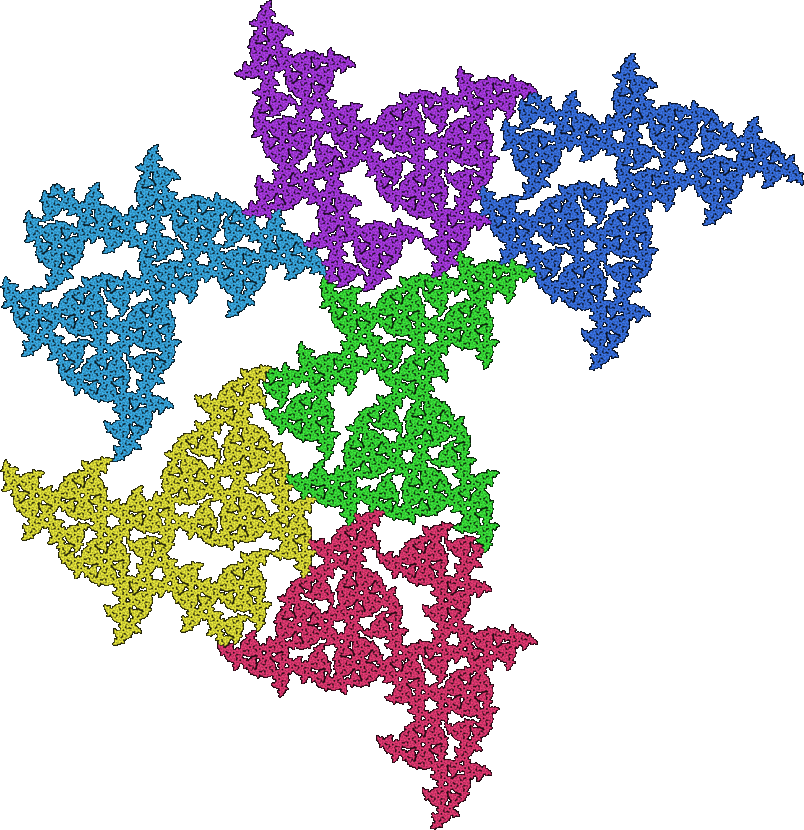} \qquad
\includegraphics[width=0.55\textwidth]{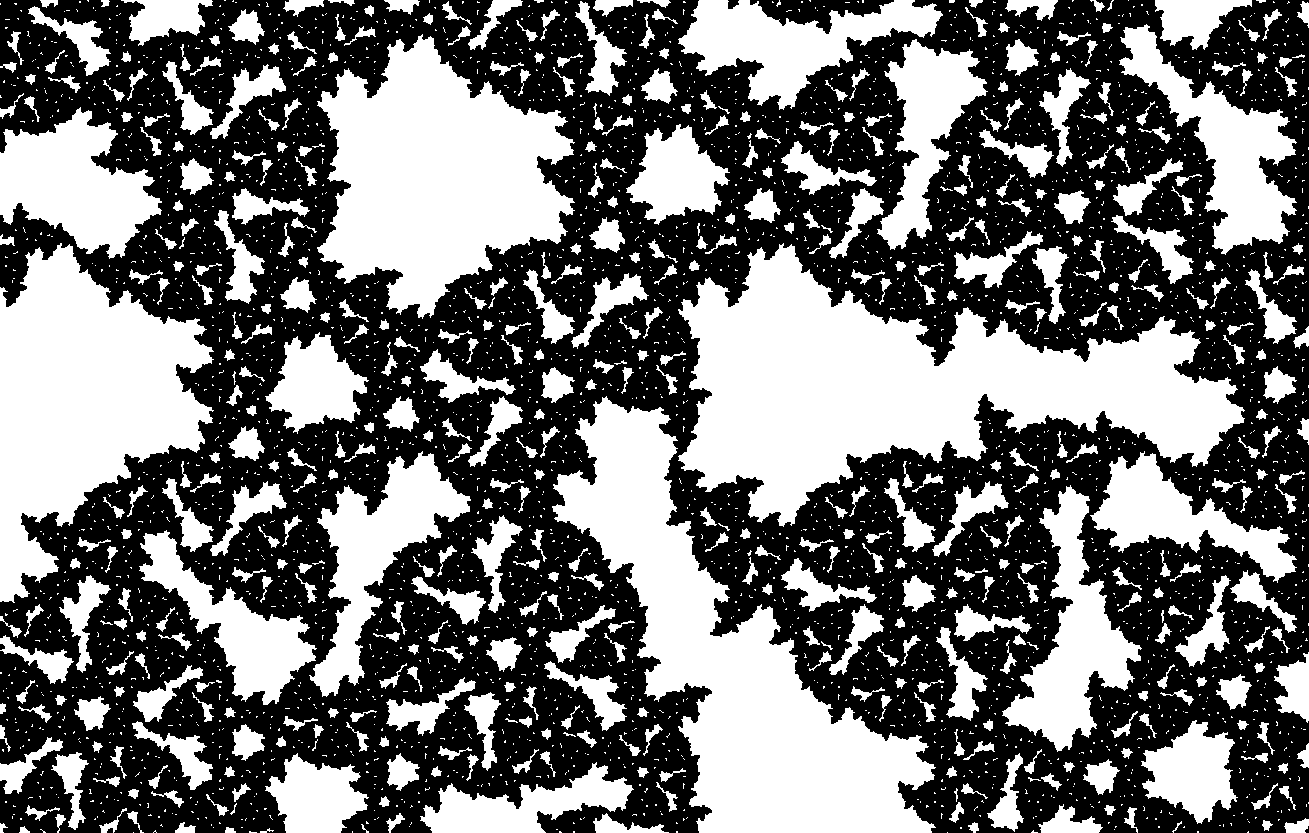}
\end{center}
\caption{A carpet generated from irrational rotation by ${\rm arctan}\,\frac{3}{13}$ together with $60^o$ rotations. The expansion has determinant 7, so for $m=6$ pieces we have dimension 1.84. There are 18 neighbor types. }\label{fi10}
\end{figure}

Figure \ref{fi10} shows the most interesting carpet with $m=6$ pieces resulting from a search with $g,s,$ and $t.$ While in most of our examples, only one piece involves the irrational rotation (relative to the other ones), here we have three  pieces with irrational rotation and three without. Similar to Figure \ref{fi6}, there seems to be no global cutpoint although local cutpoints are apparent. The close-up indicates both the symmetry of order 6 and the non-crystallographic character.  

We still have to prove that $t$ is an irrational rotation. We use the fact that $u=\frac{13}{14}, v=\frac{3}{14}$ is a rational solution of $u^2+3v^2 =1,$ and the following statement.

\begin{Proposition}\label{P03}
Let $s$ be the counterclockwise $60^o$ rotation, and let $t=as+c$ with rational $a,c.$ Let $u=c+a/2, v=a/2.$ Then $t$ is a rotation if and only if $u^2+3v^2=1.$  The rotation angle is irrational or a multiple of $60^o.$ 
\end{Proposition}

\emph{Proof. } Since $s$ has the matrix $\frac12{1 -\sqrt{3}\choose\sqrt{3}\ 1}$ with respect to the standard base, $t$ has the matrix 
${u -\sqrt{3}v\choose v\sqrt{3}\ u}.$  For a rotation, the determinant is 1. Suppose the rotation angle is rational and not a multiple of $60^o.$ Then the corresponding root of unity is the root of a cyclotomic polynomial which is irreducible over the field $\QQ (\sqrt{-3}).$ This contradicts the assumption that $t$ is a rational linear combination of $s$ and $1.$
\hfill $\Box$ \vspace{2ex}

There are a number of similar patterns in this family. As reflection for the base $B$ we can add the exchange matrix $0\ 1\choose 1\ 0$ which exchanges $b_1$ with $b_2.$ This will create a large family of carpets, related to the Gosper flake in Figure \ref{fi1}. Some of them have 100 neighbor types, but most have smaller complexity.

\section{Classification of fractal patterns}\label{class}
\subparagraph{Principles.} 
Now we shall try to get some order into our zoo of fractal examples. As in biology and crystallography, we have to divide them into species and families. This will be done in three steps.
\begin{enumerate}
\item We fix the lattice, which corresponds to the species. We can choose the square lattice with the basic $90^o$ rotation, or the hexagonal lattice with its $60^o$ rotation. Other choices are given below. Algebraically, the lattice is induced by a \emph{number field} which specifies algebraic numbers like $i$ or $\sqrt{3}$ which we need beside rational numbers. This number field comes with a symmetry group consisting of rotations defined by the number of modulus one. There can be infinitely many rotations with infinite order.
\item We choose an expansion map $g(x)=Mx$ or $g(z)=\lambda z.$ The number $\lambda$ is taken as an algebraic integer of norm greater one in our number field, sometimes as an algebraic rational. Then  $\sqrt{\det M}$ or $|\lambda|$
is the contraction factor $r$ of the IFS. Since we require OSC, the number of maps in the plane is bounded by $mr^2\le 1,$ or $m\le \det M.$  
\item The choice of lattice and expansion, together with a finite selection of rotations $s,$ determines the family of fractals. In the last step, concrete instances of this family are produced by choosing particular $s_k,v_k$ for $k=1,...,m$ such that the OSC is fulfilled.
\end{enumerate}

It makes sense to fix only a maximum number of maps and let $m$ vary within the family. Often the IFS with different numbers of maps and their neighbor graphs are related. Moreover, the actual symmetry class of an example can only be determined after the $s_k$ are chosen. Usually we are not using all symmetries which are offered.  Then we follow the convention used in classifying crystallographic patterns \cite{GS}. The symmetry type of an example is  determined by the group of maps generated from the actual neighbor maps.

\subparagraph{Examples.} 
In Figure \ref{fi1}, the Sierpi\'nski triangle was drawn in a symmetric way. When we make use of this symmetry, or include $120^o$ rotations in our IFS, we are in the hexagonal lattice, or the number field $\QQ (\sqrt{-3}),$ and have a symmetry group of three rotations. However, usually we  generate a symmetric fractal in its most simple form. That is, $f_k(x)=(x+c_k)/2$ where $c_k$ is the fixed point of $f_k.$ Then the symmetry group is trivial, $s_k=id$ for all $k,$ and all choices of non-collinear points $c_k$ are isomorphic. This symmetry type consists of a single isomorphy class for $m=3.$ For $m=2$ we would add intervals. We work in $\QQ^2,$ no extension of the rational numbers is required. 
When we go to $m=4,$ we would obtain the parallelogram, the triangle, and a lot of more fragmented tiles \cite{Ba5,M}.

When we stay with $m=3$ and allow $s_k(x)=\pm x,$ we obtain three other modified Sierpi\'nski gaskets which are well known. Again, the choice of the $c_k$ is irrelevant. We stay in $\QQ^2,$ and the family remains very small. 

In this paper, we shall always include $s(x)=-x$ in the group of symmetries, since we do not want to distinguish too many symmetry types.

In \cite{Sierrel} we studied the Sierpi\'nski triangle in a square lattice and its modifications obtained by chosing  rotations $s(z)=i^nz$ with $n\in\{ 0,1,2,3\} .$ This is a huge family. We used the square lattice $\QQ(i),$ and the symmetry group of most examples was the rotational group of the square. For $m=4,$ we get tiles like square, right-angled isosceles triangle and the aperiodic chair \cite{GS}. Of course, one gets an even larger family if reflections are added as symmetries. 

The tile on the right of Figure \ref{fi1} obviously has $120^o$ and $60^o$ rotations as neighbor maps, so we work in the hexagonal lattice, or $\QQ (\sqrt{-3}).$
The symmetry group is a rotation group, even when we define the $s_k$ with  reflections in order to keep the positive determinant of $g.$
The family is large. For $m=7$ it contains plenty of tiles \cite{M}. One example
for $m=6$ is Figure \ref{fi2}.

The new constructions in this paper have infinite symmetry groups. In Figures \ref{fi6} and \ref{fi3} the group is cyclic, with a single generator. A second generator is chosen as a reflection in Figure \ref{fi7}, and as a rotation in  Figures \ref{fi8} and \ref{fi9}. The number field is $\QQ(i),$ the lattice is formed by the Gaussian integers. For Figure \ref{fi10}, with an irrational rotation plus $60^o$ rotation, we have the hexagonal lattice of Eisenstein integers.

A large part of the literature on fractal tilings \cite{BGG,DJN12,DLN18,Ke96,KS,LL13,ScT,T89} is concerned with self-affine tiles where neighbor maps are translations and symmetries play no part. However, there are also crystallographic fractal tiles where we have a crystallographic group which acts on the tiling \cite{GG,Lo,LZ17}. Only very few fractal tiles have infinite symmetry groups \cite{Ra94,quaquaversal,Fr08,BMT}. This can be explained by a result on algebraic expansion constants which goes back to Thurston \cite{Ke96,KS,T89}.

In the study of Ngai, Sirvent, Veerman and Wang \cite{NSVW} on fractal tiles with $m=2$ pieces, symmetries play an important r\^ole. They consider rational rotations and no reflection.  They show that there are only six cases: twindragon,   L\'evy curve, Heighway dragon, triangle, rectangle and tame twindragon. The first four cases are realized in $\QQ(i)$ with expansion constant $\lambda=1+i.$ 
The hexagonal lattice does not contain a lattice point with norm $m=2,$ so it can be used only for tiles with 3 pieces.  This means that rectangle and tame twindragon are based on other lattices which we consider below.

\section{Lattices from polynomials}\label{poly}
\subparagraph{Parameter-free description of mappings.}
In this paper, a lattice was defined by a base $B=\{ b_1, b_2\}$ for which matrices $M, s_k$ have rational entries, and vectors $v_k$ have integer coordinates. How do we find $B$? Let us first reformulate condition \eqref{intas} for a rotational symmetry and expansion. 

\begin{Proposition}\label{P1}
In the plane, let $s$ be a rotation, and $g$ an expanding similarity map with positive determinant. Assume that $s$ and $g$ have a rational matrix with respect to the same base $B.$
\begin{enumerate}
\item[(i)] The characteristic polynomial of $s$ is  $p_s(z)=z^2+az+1,$ where $a$ is  rational and $|a|\le 2.$
\item[(ii)] There are rational numbers $b,c$ such that $g(x)=bs(x)+cx$ for all $x.$ 
\end{enumerate}
\end{Proposition}

\emph{Proof. }
The determinant of a rotation is 1, and the trace is rational by assumption. The eigenvalues of $s$ are $-a/2 \pm\sqrt{a^2/4 -1}.$ For $|a|>2,$ the eigenvalues are different reals so that $s$ cannot be a rotation.

For \emph{(ii)} we consider similitudes with positive determinant in standard coordinates. They have matrices of the type  ${u -v\choose v\ \ u}$ and thus form a two-dimensional vector space. The identity map $1$ together with $s$ forms a base of this space. So the
map $g$ has the form $g=b\cdot s+c\cdot 1$ with real coefficients $b,c.$ This equation is true for standard coordinates as well as for coordinates with respect to any other base $B.$
Taking a base for which both maps have a rational matrix, we see that $b$ and $c$ must be rational, by studying first values off the diagonal and then on the diagonal.
\hfill $\Box$ \vspace{2ex}

The proposition says that any triple $a,b,c$ of rational parameters defines a family of fractals. Not all triples lead to interesting examples. Our approach first selects a basic symmetry $s$ and then an expansion $g$ which fits the symmetry $s.$ One can also first define the expansion and then add symmetries in the form $s=bg+c.$ However, ordinary integer expansions, like $g(x)=2x,$ will fit any symmetry group.
Moreover, the coefficients $b,c$ for $g=bs+c$ are often integers while for $s=bg+c$ they are always rational. This led us to choose the symmetry type first. Together with a given rotation $s,$ all powers $s^q$ or $s^{-q}$ will be admitted as symmetry maps $s_k$ in \eqref{ghut}. In reality, the integer $q$ will of course be small. 

\subparagraph{Companion matrix and exchange matrix.}
The definition of a fractal family by one characteristic polynomial $p_s(z)$ and a polynomial equation $g=\sum_{k=0}^n b_ks^k$ was implemented in IFStile for a higher-dimensional setting.  The user has only to think about the master equations, while all matrix calculations are done by the computer.   
In the case when the polynomial $p_s(z)$ is irreducible, the canonical base $B$ is given by the successive images $b_k=s^{k-1}(b_1),\ k=2,...,n$ of the first base vector $b_1.$ In other words, we work with the \emph{companion matrix} of $s.$ An advantage is that the transform to standard coordinates can be done by fast numerical procedures.

Here we work with quadratic polynomials $p_s(z)=z^2+az+1$ which are irreducible for $|a|<2$ while for $a=\pm 2$ there are no interesting symmetries. We define the companion base as $b_1={1\choose 0}, b_2=s(b_1).$ For this base, $s$ has the companion matrix  $M_s={0 -1\choose 1 -a}$. Because of the master equation  $g= bs+c,$ the matrix of $g$ in this base is  $M_g=bM_s+cI$ where $I$ denotes the unity matrix. Thus  
\begin{equation}\det M=b^2+c^2-abc\quad \mbox{ and } \quad {\rm trace\,} M =2c-ab\ .
\label{detM}\end{equation}

So far we discussed only rotations and expansions with positive determinant. What about reflections? It turns out that we need only add a single reflection $r$ to our symmetry group.  Further reflections are produced from composition of $r$ with rotations. When working with the companion base, a canonical choice for $r$ is the exchange matrix $M_r={0 \ 1\choose 1 \ 0}.$ In other words, $r(b_1)=b_2$ and $r(b_2)=b_1.$ Other definitions of a basic reflection $r$ are possible. For $\QQ(i),$ ordinary conjugation is the canonical choice.

With a reflection $r,$ we can also describe IFS with expansion maps with negative determinant. Note that the definition of an attractor by an IFS $g, h_1,..,h_m$ by \eqref{ghut} is by no means unique. For any isometry $t,$ we can pass to the IFS  $tg, th_1,...,th_m$ which has the same attractor $A.$ Moreover, we can change the coordinate system and permute the maps of the IFS. The recognition of isomorphic IFS representations is still a problem of the IFStile software, in particular for symmetric examples as square and cube. See \cite[Section VI]{Sierrel} for a brief discussion of this problem.  We try to choose the simplest IFS for a given attractor, but the computer will not always do.

Finally, we mention that in the two-dimensional case we have little problems with commutativity of matrix multiplication: $s$ commutes with $g,$ and different rotations commute with each other. For a reflection $r$ and a rotation $s$ we have $sr=rs^{-1}.$

\begin{figure}[h!] 
\begin{center}
\includegraphics[width=0.34\textwidth]{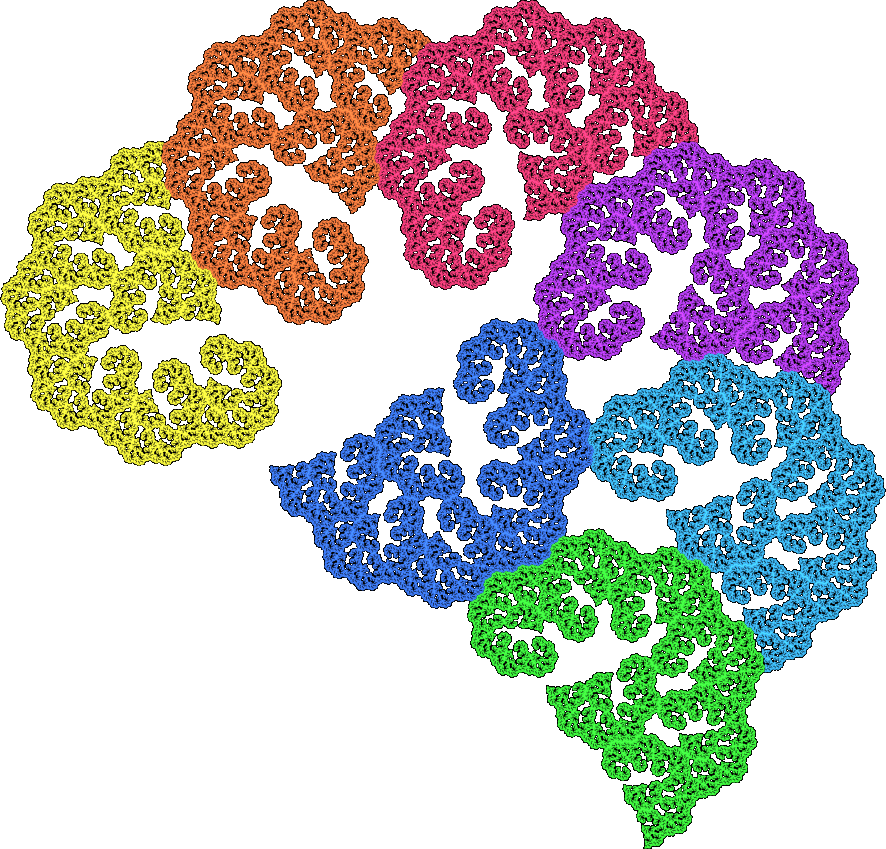} \qquad
\includegraphics[width=0.54\textwidth]{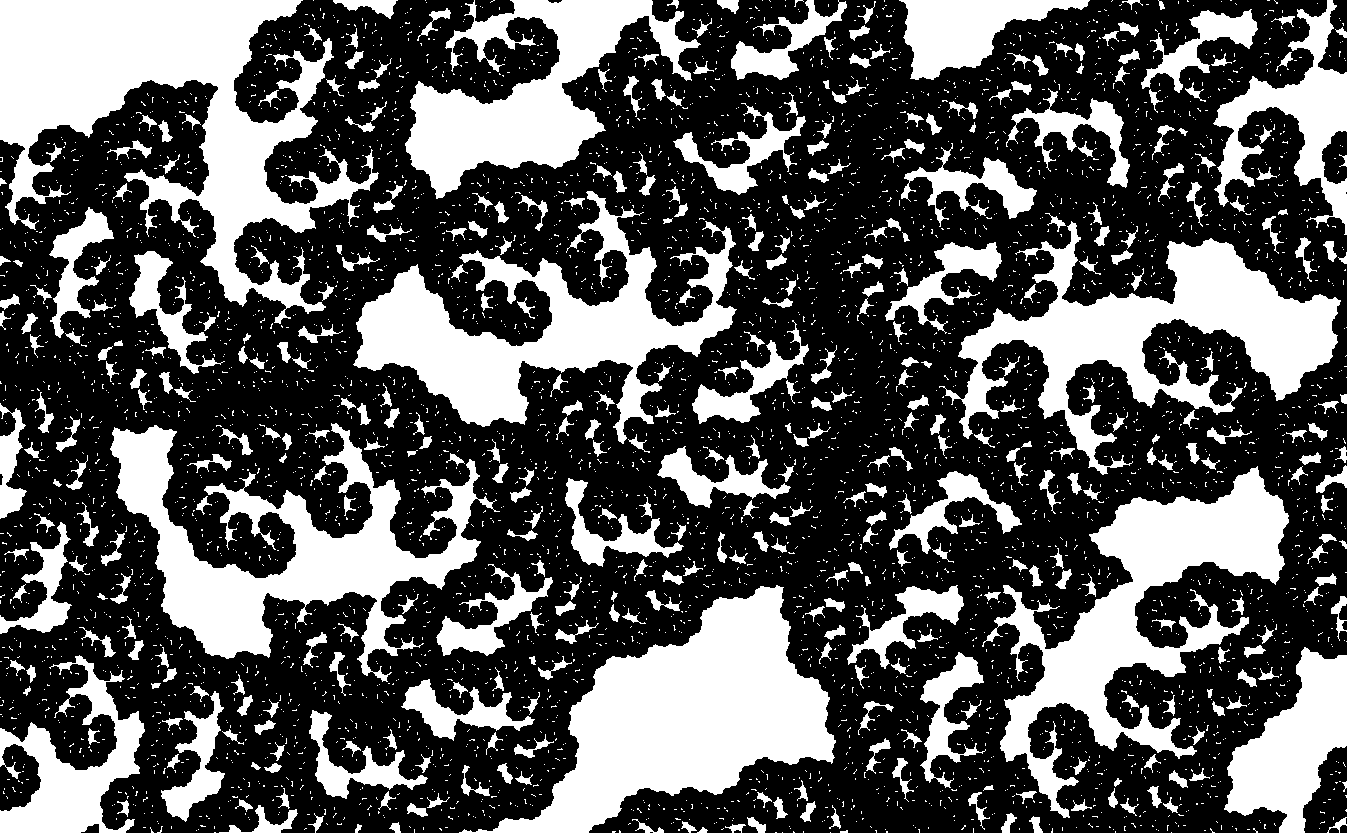}
\end{center}
\caption{A carpet generated from three irrational rotations with $a=\frac{3}{2}, \frac14, \frac98$ and an expansion $g$ of determinant 8. It has dimension 1.87 and 19 neighbor types. The number field is $\QQ(\sqrt{-7}),$ as for the tame twindragon. }\label{fi11}
\end{figure}

\section{Constructions in quadratic number fields}\label{quad}
\subparagraph{Quadratic number fields.}
Every irreducible polynomial $p$ with rational coefficients gives rise to a field extension of $\QQ .$
In particular, $p(z)=z^2+d$ with a positive square-free integer $d$ generates the field
\begin{equation}\QQ(\sqrt{-d})= \{ c_1+c_2\cdot i\sqrt{d}\ |\ c_1,c_2\in \QQ \}\ .   
\label{Qd}\end{equation}
The integer $d$ is square-free if it is not divisible by a square of an integer $n>1.$ For other integers $d$ such an $n$ would be included in the coefficient $c_2.$  
Any quadratic polynomial with rational coefficients and complex roots will generate one of these fields, as can be easily checked by solving the quadratic equation. 

Here we are interested only in polynomials $p_s(z)=z^2+az+1$ of rotations with a positive rational number $a<2.$ The number $-a$ gives the same quadratic field, and is included in any IFS with $s$ since the map $-x$ is always adjoined to the symmetry group.  We take positive integers $u<w$ with $a=2u/w.$ Then
\[ z_1=\frac{u+i\sqrt{w^2-u^2}}{w}= \frac{u+v\sqrt{-d}}{w} \]
where $w^2-u^2=dv^2.$ The number $d$ is taken as square-free part of $w^2-u^2.$ Table \ref{ad} shows all parameters $a$ with denominator $w\le 9$ which correspond to $d\le 35.$  We have characterized the rotations which belong to quadratic number fields:

\begin{Proposition}\label{P2}
Consider the polynomial $p_s(z)=z^2+az+1$ with rational $a\in (0,2).$  Let $a=2u/w$ where $u,w$ are positive integers. 
\begin{enumerate}
\item[(i)] $p_s$ generates the quadratic number field $\QQ(\sqrt{-d})$ where $d$ is the square-free part of $w^2-u^2.$ 
\item[(ii)] If $a$ is different from 0 and 1, the rotation angle $\alpha$ is irrational and fulfils $\cos\alpha = -a/2.$
\item[(iii)] The parameter \ $a$ \ corresponds to the integer triple $(u,v,w)$ which solves the equation  $u^2+ dv^2 =w^2.$ For every $d,$ these generalized Pythagorean triples are generated by Euclid's formula
\[ u=n^2-dm^2\, ,\quad v= 2mn\, , \quad w=n^2+dm^2 \quad\mbox{ with integers }m,n \mbox{ such that } n>m\sqrt{d}.\]
\end{enumerate} 
\end{Proposition}

\emph{Proof. }
The first assertion was proved above, the second follows from the fact that the cyclotomic polynomials which generate rational rotations have integer coefficients.
The rotation angle is defined by $z_1=\cos\alpha +i\sin\alpha.$ The quadratic Diophantine equation was derived above. It is equivalent to $(w-u)(w+u)=v^2d.$
Now take relatively prime integers $m,n$ with
\[ \frac{w+u}{v}= \frac{n}{m} = \frac{vd}{w-u}\quad\mbox{ and thus }\quad 
 \frac{w-u}{v}= \frac{md}{n}\ .\]
Sum and difference of both equations yields $ \frac{w}{v}= \frac{n^2+dm^2}{2mn}$ and $\frac{u}{v}= \frac{n^2-dm^2}{2mn},$ respectively. Since $m,n$ are relatively prime and $d$ is square-free, this implies (iii).
\hfill $\Box$ \vspace{2ex}

\begin{table}[h]
\centering\Large
\begin{tabular}{|c|c||c|c||c|c|c|} \hline \rule[-3mm]{0mm}{8mm}
$d \, $&$a$&$d$&$a$&$d$&$a$\\ \hline \rule[-5mm]{0mm}{14mm}
\ 1 \ &$\displaystyle 0,\ \frac{6}{5},\ \frac{8}{5}$&\ 10 \ &$\displaystyle\frac{6}{7}$&\ 19 \ &$\displaystyle\frac{9}{5}$\\ \hline \rule[-5mm]{0mm}{14mm}
\ 2 \ &$\displaystyle \frac{2}{3}, \ \frac{14}{9}$&11&$\displaystyle\frac{5}{3}, \ \frac{1}{5}, \ \frac{7}{9}$&21&$\displaystyle\frac{4}{5}$\\ \hline
\rule[-5mm]{0mm}{14mm}
\ 3 \ &$\displaystyle 1, \frac{2}{7}, \frac{11}{7}, \frac{13}{7}$&13&$\displaystyle \frac{12}{7}$&23&$\displaystyle\frac{11}{6}, \ \frac{7}{8}$\\ \hline
\rule[-5mm]{0mm}{14mm}
\ 5 \ &$\displaystyle \frac{4}{3}, \ \frac{4}{7}, \ \frac{2}{9}$&14&$\displaystyle \frac{10}{9}$&31&$\displaystyle \frac{15}{8}$\\ \hline
\rule[-5mm]{0mm}{14mm}
\ 6 \ &$\displaystyle  \ \frac{2}{5}, \ \frac{10}{7}$&15&$\displaystyle\frac12 , \ \frac{7}{4}, \ \frac{11}{8}$&33&$\displaystyle\frac{8}{7}$\\ \hline
\rule[-5mm]{0mm}{14mm}
\ 7 \ &$\displaystyle \frac{3}{2}, \ \frac{7}{4}, \ \frac{9}{8}$&17&$\displaystyle \frac{16}{9}$&35&$\displaystyle\frac{1}{3}, \ \frac{17}{9}$\\ \hline
\end{tabular}
\caption{The first positive square-free integers $d,$ and corresponding fractions $a$ with denominator up to 9 which generate the field $\QQ(\sqrt{-d}).$ }\label{ad}
\end{table}

\subparagraph{Results of the computer experiments.}
We studied rational numbers $a$ with small denominator and $|a|<2,$ and the rotation with polynomial $p_s(z)=z^2+az+1.$ The expansion was defined as $g=bs+c$ with rational numbers. The best results were obtained when $g$ represents an algebraic integer in $\QQ(\sqrt{-d}).$ That is, determinant and trace of $M$ in \eqref{detM} must be integers. Recall that an element of the field \eqref{Qd} is an algebraic integer if either $c_1,c_2$ are integers, or $d=3 {\rm mod }4$ and $c_1=d_1/2,\ c_2=d_2/2$ with odd integers $d_1,d_2.$

For $a=\frac32 ,$ we get the field $\QQ(\sqrt{-7})$ which corresponds to the tame twindragon \cite{Ba5,NSVW}. The twindragon is obtained for $g=2s+1$ but does not use the symmetry $s.$ With $g=3s+1,$ we got an example with $m=4$ which is extremely similar to Figure \ref{fi3}. It has 6 neighbor types and fulfils $\beta =\alpha/4$ as in Section \ref{ngapp}. This is remarkable since this $g$ has determinant 11/2 and trace 5/2 and thus is not an algebraic integer.

For the expansion $g=4s+1$ with determinant 11, we got many nice carpets with $m=9$ and small complexity. Adding reflections, we got an even greater variety of shapes.  One example had dimension $\alpha=1.83$ and almost the same boundary dimension $\beta=1.76,$ similar to the L\'evy curve \cite{DKV,SW}.  Another interesting expansion was $g=2s-2$ with determinant 14. 

Figure \ref{fi11} shows a very interesting example where all three different irrational rotations from Table \ref{ad} are involved. Here we took $g=2s-1$ with determinant 8, and we have $m=7$ pieces. The Hausdorff dimension is 1.87, and there are no cutpoints whatsoever. Only two pieces are parallel and two pieces related by point reflection. All other pairs of pieces differ by different irrational rotations. From all our IFSs, the isotropic character is most obvious here.
Nevertheless we have only 19 neighbors. There were two variations of this example, one with 16 neighbors and cutpoints, another one with 21 neighbors and still better geometry. On the whole $a=3/2$ with $\QQ(\sqrt{-7})$ is a rich source of irrational examples.

\begin{figure}[h!] 
\begin{center}
\includegraphics[width=0.34\textwidth]{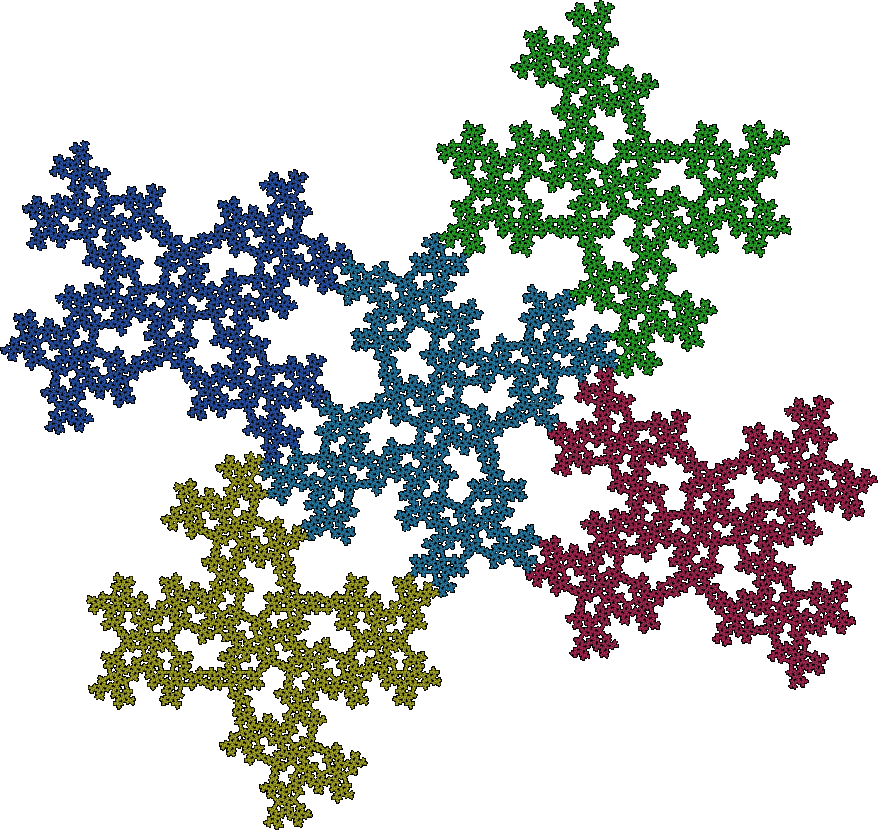} \qquad
\includegraphics[width=0.51\textwidth]{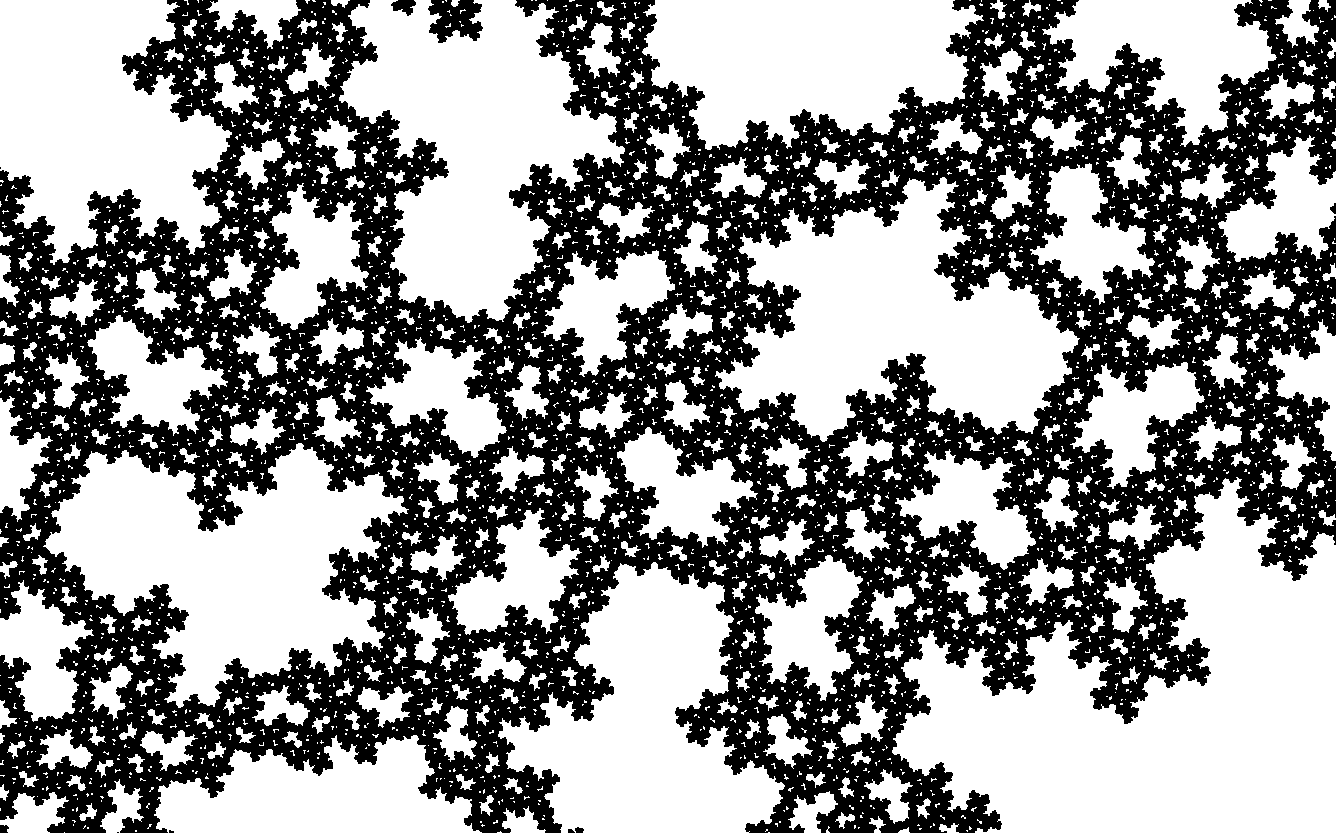}
\end{center}
\caption{One of the simplest carpets for $a=1/2,$ which corresponds to $\QQ(\sqrt{-15}).$  }\label{fi12}
\end{figure}

For $a=1/2,$ which corresponds to $\QQ(\sqrt{-15}),$ we got also rich families. The simplest expansion is $g=2s-1$ with determinant 6. Of course we did not find tiles, but abundant examples with $m=5.$ Similar to crystallographic IFS like Figure \ref{fi2}, the complexity can be large, and the pieces of $A$ may cross each other which indicates that the OSC can only be fulfilled with complicated open sets \cite{Sierrel}. An example with only 12 neighbors is shown in Figure \ref{fi12}. The rotation with $a=7/4$ leads to the same number field and similar examples.

The rectangle as a self-similar set with two pieces, used as standard size for writing paper, belongs to the field $\QQ(\sqrt{-2}).$ An appropriate rotation is given by $a=2/3.$ The expansion $g=2s-1$ provides one example with $m=4$ similar to Figure \ref{fi3}, but with smaller dimension 1.5 and larger number 12 of neighbors. The reason might be that $g$ with determinant 19/3 is not an algebraic integer. Another expansion is $g=s-3$ with determinant 12. We obtained connected examples with OSC only up to $m=8.$ Still these examples have global cutpoints, like Figure \ref{fi3}, even if reflections are used. A third choice was $g=3s+1$ with determinant 8. The search gave a single example with $m=7$ and Cantor sets as boundaries. It is exceptional since it has only 6 neighbor types. But like Figure \ref{fi3}, it has global cutpoints and involves a tree structure. Altogether, the field 
$\QQ(\sqrt{-2})$ seems not to lead to nice carpets for $m\le 12.$

We performed similar experiments for $d=5,11,19,23,35,91$ and got similar results as for $d=2,$ at least as good as Figure \ref{fi3}. For $d=6$ there were less examples.

\subparagraph{Summary.}
Including irrational rotations into IFS is not easy. No tiles could be constructed this way. Nevertheless, we found connected self-similar sets with pieces intersecting in uncountable sets in all quadratic fields which we studied. Various examples have the topology of the Sierpi\'nski carpet, almost the same Hausdorff dimension, and a much more interesting isotropic geometry. For the Sierpi\'nski triangle, with finite intersections of pieces and dimension between 1.5 and 1.6, there are lots of variations in all quadratic fields considered.
In general, the complexity of irrational examples is smaller than that of IFS with crystallographic data. Some of the examples with small number of neighbors seem to be unique and worth of further study.

\subparagraph{Some open problems.} 
This is a beginning. We discussed self-similar sets in the plane with equal  contraction factors, between 4 and 14 pieces, and data from quadratic number fields. Graph-directed constructions and projection schemes for arbitrary number fields are more exciting. For three-dimensional self-similar sets a general framework does not yet exist, despite much work, including \cite{Ba10,CT16,quaquaversal,DLN18,Encyc,KS,M,TZ20}. 

In Section \ref{intro} we claimed that our examples do not contain line segments. How can this be proved?  If there is a line segment $L$ in $A,$ it will intersect a small piece $A_w$ in an interior subsegment $L_0,$ and meet two subpieces  $A_{wj}$ and $A_{wk}.$ Then $L'=f_w^{-1}(L_0\cap A_w)$ is a segment which connects two boundary points of $A$ in different pieces $A_j, A_k.$ In all our figures, visual inspection shows that such $L'$ cannot exist. It would be better to prove this by a characterization of all IFS which lead to segments in $A.$

It should also be possible to characterize attractors with cutpoints and give a fast algorithm for detecting Cantor set attractors, which may have very large neighbor graphs \cite{Sierrel}. The holes of all irrational examples obviously look more complicated than the holes of Sierpi\'nski's spaces. However, calculation of the Hausdorff dimension of the topological boundary of such holes is still an open problem.

\vspace{3ex}

\noindent
Christoph Bandt\\  Institute of Mathematics, University of Greifswald, 17487 Greifswald, Germany. \\ 
\url{bandt@uni-greifswald.de}\vspace{1ex}\\
Dmitry Mekhontsev\\Sobolev Institute of Mathematics, 
4 Acad. Koptyug avenue, 630090 Novosibirsk Russia\\
\url{mekhontsev@gmail.com}\\


\begin{thebibliography}{10}

\bibitem{AL}
S.~Akiyama and B.~Loridant.
\newblock Boundary parametrization of self-affine sets.
\newblock {\em J. Math. Soc. Japan}, 63(2):525--579, 2011.

\bibitem{Baake2013}
M.~Baake and U.~Grimm.
\newblock {\em Aperiodic Order, Vol. 1: A mathematical invitation}.
\newblock Cambridge University Press, Cambridge, 2013.

\bibitem{Ba5}
C.~Bandt.
\newblock Self-similar sets 5. {Integer} matrices and fractal tilings of
  {${\mathbb R}^n$}.
\newblock {\em Proc. Amer. Math. Soc.}, 112:549--562, 1991.

\bibitem{Ba10}
C.~Bandt.
\newblock Combinatorial topology of three-dimensional self-affine tiles.
\newblock 2010.
\newblock arXiv:1002.0710.

\bibitem{BGG}
C.~Bandt and G.~Gelbrich.
\newblock Classification of self-affine lattice tilings.
\newblock {\em J. London Math. Soc.}, 50:581--593, 1994.

\bibitem{BG}
C.~Bandt and S.~Graf.
\newblock Self-similar sets 7. {A} characterization of self-similar fractals
  with positive {Hausdorff} measure.
\newblock {\em Proc. Amer. Math. Soc.}, 114:995--1001, 1992.

\bibitem{BHR}
C.~Bandt, N.V. Hung, and H.~Rao.
\newblock On the open set condition for self-similar fractals.
\newblock {\em Proc. Amer. Math. Soc.}, 134:1369--1374, 2005.

\bibitem{Sierrel}
C.~Bandt and D.~Mekhontsev.
\newblock Elementary fractal geometry. new relatives of the sierpinski gasket.
\newblock {\em Chaos: An Interdisciplinary Journal of Nonlinear Science},
  28(6):063104, 2018.

\bibitem{BM19}
C.~Bandt and D.~Mekhontsev.
\newblock Computer geometry: Rep-tiles with a hole.
\newblock {\em Mathematical Intelligencer}, 42(1):1--5, 2020.

\bibitem{BMT}
C.~Bandt, D.~Mekhontsev, and A.~Tetenov.
\newblock A single fractal pinwheel tile.
\newblock {\em Proc. Amer. Math. Soc.}, 146:1271--1285, 2018.

\bibitem{BM09}
C.~Bandt and M.~Mesing.
\newblock Self-affine fractals of finite type.
\newblock In {\em Convex and fractal geometry}, volume~84 of {\em Banach Center
  Publ.}, pages 131--148. Polish Acad. Sci. Inst. Math., Warsaw, 2009.

\bibitem{Bar}
M.~F. Barnsley.
\newblock {\em Fractals everywhere}.
\newblock Academic Press, 2nd edition, 1993.

\bibitem{BP}
C.J. Bishop and Y.~Peres.
\newblock {\em Fractal sets in probability and analysis}.
\newblock Cambridge University Press, Cambridge, 2017.

\bibitem{CT16}
G.~R. Conner and J.~M. Thuswaldner.
\newblock Self-affine manifolds.
\newblock {\em Advances Math.}, 289:725--783, 2016.

\bibitem{quaquaversal}
J.H. Conway and C.~Radin.
\newblock Quaquaversal tilings and rotations.
\newblock {\em Invent. Math.}, 132:179--188, 1998.

\bibitem{DE11}
M~Das and GA~Edgar.
\newblock Finite type, open set conditions and weak separation conditions.
\newblock {\em Nonlinearity}, 24(9):2489, 2011.

\bibitem{DJN12}
D.-W. Deng, T.~Jiang, and S.-M. Ngai.
\newblock Structure of planar integral self-affine tilings.
\newblock {\em Mathematische Nachrichten}, 285(4):447--475, 2012.

\bibitem{DLN18}
G.~Deng, C.~Liu, and S.-M. Ngai.
\newblock Topological properties of a class of self-affine tiles in {${\mathbb
  R}^3$}.
\newblock {\em Transactions of the American Mathematical Society},
  370(2):1321--1350, 2018.

\bibitem{DKV}
P.~Duvall, J.~Keesling, and A.~Vince.
\newblock The {Hausdorff} dimension of the boundary of a self-similar tile.
\newblock {\em J. London Math. Soc.}, 61:748--760, 2000.

\bibitem{Fal}
K.~J. Falconer.
\newblock {\em Fractal geometry: mathematical foundations and applications}.
\newblock J. Wiley \& sons, 3 edition, 2014.

\bibitem{FFJ15}
K.~Falconer, J.~Fraser, and X.~Jin.
\newblock Sixty years of fractal projections.
\newblock In {\em Fractal geometry and stochastics V}, pages 3--25. Springer,
  2015.

\bibitem{Fr08}
D.~Frettl\"{o}h.
\newblock Substitution tilings with statistical circular symmetry.
\newblock {\em European J. Comb.}, 29:1881--1893, 2008.

\bibitem{Encyc}
D.~Frettl\"{o}h.
\newblock Tilings {Encyclopedia}, 2018.
\newblock \url{tilings.math.uni-bielefeld.de}.

\bibitem{GG}
G.~Gelbrich.
\newblock Crystallographic reptiles.
\newblock {\em Geometria Dedicata}, 51:235--256, 1994.

\bibitem{Gi86}
W.J. Gilbert.
\newblock The fractal dimension of sets derived from complex bases.
\newblock {\em Canad. Math. Bull.}, 29(4):495--500, 1986.

\bibitem{GS}
B.~Gr\"{u}nbaum and G.C. Shephard.
\newblock {\em Patterns and Tilings}.
\newblock Freeman, New York, 1987.

\bibitem{HLR}
X.-G. He, K.-S. Lau, and H.~Rao.
\newblock Self affine sets and graph-directed systems.
\newblock {\em Constr. Approx.}, 19:373--397, 2003.

\bibitem{HS12}
M.~Hochman and P.~Shmerkin.
\newblock Local entropy averages and projections of fractal measures.
\newblock {\em Annals of Mathematics}, pages 1001--1059, 2012.

\bibitem{Ke96}
R.~Kenyon.
\newblock The construction of self-similar tilings.
\newblock {\em Geom. Funct. Anal.}, 6:471--488, 1996.

\bibitem{KS}
R.~Kenyon and B.~Solomyak.
\newblock On the characterization of expansion maps for self-affine tilings.
\newblock {\em Discrete Comput. Geom.}, 43:577--593, 2010.

\bibitem{Kig}
J.~Kigami.
\newblock {\em Analysis on fractals}, volume 143.
\newblock Cambridge University Press, 2001.

\bibitem{LN07}
K.-S. Lau and S.-M. Ngai.
\newblock A generalized finite type condition for iterated function systems.
\newblock {\em Advances in Mathematics}, 208(2):647--671, 2007.

\bibitem{LL13}
K.-S. Leung and J.J. Luo.
\newblock Boundaries of disk-like self-affine tiles.
\newblock {\em Discrete Comput. Geom.}, 50(1):194--218, 2013.

\bibitem{Lo}
B.~Loridant.
\newblock Crystallographic number systems.
\newblock {\em Monatsh. Math.}, 167:511--529, 2012.

\bibitem{LZ17}
B.~Loridant and S.-Q. Zhang.
\newblock Topology of a class of p2-crystallographic replication tiles.
\newblock {\em Indagationes Mathematicae}, 28(4):805--823, 2017.

\bibitem{Man}
B.B. Mandelbrot.
\newblock {\em The fractal geometry of nature}.
\newblock Freeman, New York, 1982.

\bibitem{MW}
R.D. Mauldin and S.C. Williams.
\newblock Hausdorff dimension in graph-directed constructions.
\newblock {\em Trans. Amer. Math. Soc.}, 309:811--829, 1988.

\bibitem{M}
D.~Mekhontsev.
\newblock {IFS} tile finder, version 2.2.
\newblock \url{https://ifstile.com}, 2019.

\bibitem{Mo}
P.A.P. Moran.
\newblock Additive functions of intervals and hausdorff measure.
\newblock {\em Math. Proc. Cambridge Phil. Soc.}, 42:15--23, 1946.

\bibitem{NSVW}
S.-M. Ngai, V.F. Sirvent, J.J.P. Veerman, and Y.~Wang.
\newblock On 2-reptiles in the plane.
\newblock {\em Geometriae Dedicata}, 82(1-3):325--344, 2000.

\bibitem{NW}
S.-M. Ngai and Y.~Wang.
\newblock Hausdorff dimension of self-similar sets with overlaps.
\newblock {\em J. London Math. Soc.}, 63:655--672, 2001.

\bibitem{Ra94}
C.~Radin.
\newblock The pinwheel tilings of the plane.
\newblock {\em Annals Math.}, 139:661--702, 1994.

\bibitem{RZ16}
H.~Rao and Y.~Zhu.
\newblock Lipschitz equivalence of fractals and finite state automaton.
\newblock {\em arXiv preprint arXiv:1609.04271}, 2016.

\bibitem{ScT}
K.~Scheicher and J.M. Thuswaldner.
\newblock Neighbors of self-affine tiles in lattice tilings.
\newblock In P.~Grabner and W.~Woess, editors, {\em Fractals in Graz 2001},
  pages 241--262. Birkh\"auser, 2003.

\bibitem{Sch}
A.~Schief.
\newblock Separation properties for self-similar sets.
\newblock {\em Proc. Amer. Math. Soc.}, 122:111--115, 1994.

\bibitem{Se95}
M.~Senechal.
\newblock {\em Quasicrystals and geometry}.
\newblock Cambridge University Press, Cambridge, 1995.

\bibitem{WS}
S.~Sheffield and W.~Werner.
\newblock Conformal loop ensembles: the markovian characterization and the
  loop-soup construction.
\newblock {\em Annals of Mathematics}, pages 1827--1917, 2012.

\bibitem{SS18}
P.~Shmerkin and V.~Suomala.
\newblock {\em Spatially independent martingales, intersections, and
  applications}, volume 251.
\newblock American Mathematical Society, 2018.

\bibitem{Shmerkin15}
P.~Shmerkin.
\newblock Projections of self-similar and related fractals: a survey of recent
  developments.
\newblock In {\em Fractal geometry and stochastics V}, pages 53--74. Springer,
  2015.

\bibitem{SoUniq}
B.~Solomyak.
\newblock Nonperiodicity implies unique composition for self-similar
  translationally finite tilings.
\newblock {\em Discrete Comput. Geom.}, 20:265--279, 1998.

\bibitem{Str}
R.S. Strichartz.
\newblock {\em Differential equations on fractals: a tutorial}.
\newblock Princeton University Press, 2006.

\bibitem{SW}
R.S. Strichartz and Y.~Wang.
\newblock Geometry of self-affine tiles 1.
\newblock {\em Indiana Univ. Math. J.}, 48:1--24, 1999.

\bibitem{T89}
W.P. Thurston.
\newblock Groups, tilings, and finite state automata.
\newblock AMS Colloquium Lectures, Boulder, CO, 1989.

\bibitem{TZ20}
J.~Thuswaldner and S.~Zhang.
\newblock On self-affine tiles whose boundary is a sphere.
\newblock {\em Transactions of the American Mathematical Society},
  373(1):491--527, 2020.

\bibitem{Zer}
M.P.W. Zerner.
\newblock Weak separation properties for self-similar sets.
\newblock {\em Proc. Amer. Math. Soc.}, 124:3529--3539, 1996.

\end{thebibliography}
\end{document}